\documentclass[preprint]{imsart}
\usepackage{xr}
\RequirePackage[OT1]{fontenc}
\RequirePackage{amsthm,amsmath,natbib}
\RequirePackage{natbib,color}
\RequirePackage[colorlinks,citecolor=blue,urlcolor=blue]{hyperref}
\RequirePackage[colorlinks]{hyperref}
\startlocaldefs
\numberwithin{equation}{section}
\theoremstyle{plain}

\endlocaldefs
\usepackage{psfrag, epsfig, graphicx}
\usepackage{amsfonts,amsmath,latexsym,amssymb}
\usepackage[active]{srcltx}
\usepackage[a4paper]{geometry}
\newtheorem{lemma}{Lemma}
\newtheorem{theorem}{Theorem}
\newtheorem{proposition}{Proposition}

\newtheorem{assumption}{Assumption}
\newtheorem{remark}{Remark}
\newtheorem{definition}{Definition}
\newtheorem{corollary}{Corollary}

\newtheorem{assumption*}{Assumption}

\renewcommand{\kappa}{\varkappa}

\newcommand{\rd}{{\rm d}}
\newcommand{\e}{\varepsilon}

\newcommand{\cA}{{\cal A}}
\newcommand{\cB}{{\cal B}}
\newcommand{\cC}{{\cal C}}
\newcommand{\cD}{{\cal D}}
\newcommand{\cE}{{\cal E}}
\newcommand{\cF}{{\cal F}}

\newcommand{\cI}{{\cal I}}
\newcommand{\cJ}{{\cal J}}

\newcommand{\cM}{{\cal M}}

\newcommand{\cR}{{\cal R}}

\newcommand{\cW}{{\cal W}}
\newcommand{\cX}{{\cal X}}

\newcommand{\Ba}{\boldsymbol{a}}

\newcommand{\Bc}{\boldsymbol{c}}

\newcommand{\Be}{\boldsymbol{e}}

\newcommand{\BL}{\boldsymbol{L}}

\newcommand{\BM}{\boldsymbol{M}}

\newcommand{\Bt}{\boldsymbol{t}}

\newcommand{\Bu}{\boldsymbol{u}}

\newcommand{\bma}{\boldsymbol{\sigma}}

\newcommand{\bze}{\zeta}

\newcommand{\bB}{\mathbb B}
\newcommand{\bC}{\mathbb C}

\newcommand{\bE}{\mathbb E}
\newcommand{\bF}{\mathbb F}

\newcommand{\bL}{{\mathbb L}}

\newcommand{\bN}{{\mathbb N}}

\newcommand{\bP}{{\mathbb P}}

\newcommand{\bR}{{\mathbb R}}

\newcommand{\bX}{{\mathbb X}}

\newcommand{\mA}{A}
\newcommand{\mB}{\mathfrak{B}}

\newcommand{\mD}{\mathfrak{D}}
\newcommand{\mF}{\mathfrak{F}}

\newcommand{\mL}{\mathfrak{L}}
\newcommand{\mP}{\mathfrak{P}}

\newcommand{\mJ}{\mathfrak{J}}

\newcommand{\mn}{\mathfrak{n}}

\newcommand{\mT}{\mathfrak{T}}
\newcommand{\mX}{\mathfrak{X}}

\newcommand{\md}{\mathfrak{d}}

\newcommand{\mR}{\mathfrak{R}}
\newcommand{\mz}{\mathfrak{z}}
\newcommand{\me}{\mathfrak{e}}
\newcommand{\ms}{s}
\newcommand{\ma}{\mu}
\newcommand{\mb}{\nu}
\newcommand{\mc}{\pi}
\newcommand{\mmp}{p}

\newcommand{\epr}{\hfill\hbox{\hskip 4pt
                \vrule width 5pt height 6pt depth 1.5pt}\vspace{0.5cm}\par}

\newcommand{\bva}{\boldsymbol{\varsigma}}

\begin{document}
\begin{frontmatter}
\title{Minimax estimation of norms of a probability density: I. Lower bounds}
\runtitle{Estimation of norms of a probability density}
\begin{aug}
\author
{
\fnms{A.} \snm{Goldenshluger}
\thanksref{t1}
\ead[label=e1]{goldensh@stat.haifa.ac.il}
}
\and
\author
{\fnms{O. V.} \snm{Lepski}
\thanksref{t2}
\ead[label=e2]{oleg.lepski@univ-amu.fr}}
\runauthor{A.~Goldenshluger and O. V.~Lepski}

\affiliation{University of Haifa\thanksmark{m1} and Aix--Marseille Universit\'e, CNRS, Centrale Marseille,
I2M\thanksmark{m2}}

\address{Department of Statistics\\
University of Haifa
\\
Mount Carmel
\\ Haifa 31905, Israel\\
\printead{e1}}


\address{Institut de Math\'ematique de Marseille\\
Aix-Marseille  Universit\'e   \\
 39, rue F. Joliot-Curie \\
13453 Marseille, France\\
\printead{e2}}
\end{aug}

\thankstext{t1}{Supported by the ISF grant No. 361/15.}
\thankstext{t2}{This work has been carried out in the framework of the Labex Archim\`ede (ANR-11-LABX-0033) and of the A*MIDEX project (ANR-11-IDEX-0001-02), funded by the "Investissements d'Avenir" French Government program managed by the French National Research Agency (ANR).}
\begin{abstract}
The paper deals with the problem of
nonparametric estimating the $\bL_p$--norm, $p\in (1,\infty)$, of a probability density
on $\bR^d$, $d\geq 1$ from independent observations.
The unknown density
is assumed to belong to a ball in the anisotropic
Nikolskii's space. We adopt the minimax approach, and derive
lower bounds on the minimax risk.
In particular,
we demonstrate that accuracy of estimation procedures
essentially depends  on whether $p$ is integer or not. Moreover, we develop
a general technique for derivation of lower bounds on the minimax risk
in the problems of estimating nonlinear functionals. The proposed
technique is applicable
for a broad class of nonlinear functionals, and it is used for derivation of
the lower bounds in the~$\bL_p$--norm estimation.
\end{abstract}


\end{frontmatter}

\section{Introduction}
 Suppose that we observe i.i.d. vectors $X_i\in\bR^d, i=1,\ldots, n,$ with
 common probability density $f$.
 Let $p>1$ be a given real number. We want to estimate the $\bL_p$-norm of $f$,
 $$
 \|f\|_p:=\bigg[\int_{\bR^d}|f(x)|^p\rd x\bigg]^{1/p},
 $$
using observations $X^{(n)}=(X_1,\ldots,X_n)$.
By estimator  we mean any $X^{(n)}$-measurable
map $\widetilde{F}:\bR^n\to\bR$, and  accuracy of an estimator $\widetilde{F}$
is measured by the quadratic risk
\[
 \cR_n[\widetilde{F}, f]:=\Big(\bE_f \big[\widetilde{F}-\|f\|_p\big]^2\Big)^{1/2},
\]
where $\bE_f$ denotes  expectation with respect to the probability measure
$\bP_f$ of  observations $X^{(n)}=(X_1,\ldots,X_n)$.
\par
We adopt the minimax approach to the outlined estimation problem.
Let $\mF$ denote the set of all probability densities defined on $\bR^d$. With any
estimator $\widetilde{F}$ and any subset $\cF$ of
$\mF$ we associate  {\em the maximal risk} of $\widetilde{F}$ on $\cF$:
$$
\cR_n\big[\widetilde{F}, \cF\big]:=\sup_{f\in\cF}\cR_n[\widetilde{F}, f].
$$
{\em The minimax risk} is
\[
 \cR_n[\cF]:=\inf_{\tilde{F}} \cR_n[\widetilde{F}, \cF],
\]
where $\inf$ is taken over all possible estimators. An estimator $\widetilde{F}_*$ is called
{\em optimal in order} or {\em rate--optimal} if
\[
 \cR_n[\widetilde{F}_*; \cF] \asymp \cR_n[\cF],\;\;\;n\to\infty.
\]
The rate at which $\cR_n[\cF]$ converges to zero as $n$ tends to infinity is referred to as
{\em the minimax rate of convergence}.
\par
The problems of minimax
nonparametric  estimation  of density functionals
have been extensively studied in the literature. The case of linear functionals is
particularly well understood: here  a
complete optimality theory under rather general assumptions has been developed
[see, e.g., \cite{Ibragimov-Khas}, \cite{Donoho-Liu}, \cite{Cai-Low-1} and
\cite{JudNem}].
As for nonlinear functionals, the situation is completely different:
even in the problem of estimating quadratic functionals of a density rate--optimal
estimators are known only for very specific functional classes.
For representative publications dealing with
estimation of quadratic and closely related
integral
functionals  of a probability density we refer to
\cite{bickel-ritov}, \cite{birge-massart}, \cite{picard}, \cite{beatrice}, \cite{gine-nickl}
and \cite{waart}.
The problems of estimating non-linear functionals were also considered in the framework of
the Gaussian white noise model; e.g., \cite{Ibragimov-et-al}, \cite{Nem90}
[see also \cite[Chapters 7 and 8]{Nem2000}], \cite{Donoho-Nuss},
\cite{cai-low}.
The contribution of this paper is closely related to the works
\cite{LNS}, \cite{cai-low-2} and \cite{han&co1},
where the problem of estimation of norms of a signal observed in Gaussian white noise
was studied.
Additional pointers to relevant work
and discussion of relations between our results and the existing literature
 are provided in
Section~\ref{sec:generic-choice}.
\par
This paper deals with
the problem of estimating  the $\bL_p$--norm
of a probability density and derives  lower bounds on asymptotics of the minimax risk
over anisotropic Nikolskii's classes $\bN_{\vec{r}, d}(\vec{\beta}, \vec{L})$
(precise definition of the functional class is given below).
In the companion paper
\cite{GL20b} we develop the corresponding rate--optimal estimators demonstrating that
the derived lower bounds are tight.
We also study how boundedness of the underlying density $f$
in some integral norm influences the estimation accuracy by considering
the minimax risk
over the functional class
$\cF= \bN_{\vec{r},d}\big(\vec{\beta},\vec{L}\big)\cap \bB_q(Q)$, where
\[
\bB_q(Q):=\big\{f:\bR^d\to\bR: \|f\|_q\leq Q\big\},\; q>1,\;Q>0.
\]
\par
The contribution of this paper is two--fold.
First, we derive lower bounds
on the minimax risk on the class $\cF$
in the problem of estimating the $\bL_p$--norm  of a probability density.
Second, we develop general machinery for derivation of lower bounds on the minimax risk
in the problems of estimating nonlinear functionals of the type
\begin{equation}\label{eq:Psi}
\Psi(f)=G\bigg(\int_{\bR^d}H\big(f(x)\big)\rd x\bigg),
\end{equation}
where $G:\bR\to\bR$ and $H:\bR_+\to\bR$  are fixed functions.
The developed  machinery is applied for  the problem of estimating $\|f\|_p$.
In order to demonstrate broad applicability
of the proposed technique
we also provide lower bounds in problems of
estimation of other nonlinear functionals of interest.
\par
The rest of this paper is structured as follows.
Section~\ref{sec:Lp-norms} presents lower bounds on the minimax risk in the problem of estimating
the $\bL_p$--norms of~$f$.
Section~\ref{sec:abstract-LBT} develops a general technique for derivation
of lower bounds in the problems of estimating nonlinear functionals of type (\ref{eq:Psi}).
The main results of these two sections
are proved in Sections~\ref{sec:proofs_Tms-norms} and~\ref{sec:proofs_Tms-general} respectively.
Appendix contains proofs of auxiliary results.
\section{Lower bounds for estimation of  the $\bL_p$--norm}
\label{sec:Lp-norms}
We start with the definition of the anisotropic Nikolskii functional classes.
Let $(\Be_1,\ldots,\Be_d)$ denote the canonical basis of $\bR^d$.
 For function $G:\bR^d\to \bR$ and
real number $u\in \bR$
{\em the first order difference operator}
with step size $u$ in direction of variable
$x_j$ is defined by
$\Delta_{u,j}G (x):=G(x+u\Be_j)-G(x),\;j=1,\ldots,d$.
By induction,
the $k$-th order difference operator with step size $u$ in direction of  $x_j$ is
\[
 \Delta_{u,j}^kG(x)= \Delta_{u,j} \Delta_{u,j}^{k-1} G(x) = \sum_{l=1}^k (-1)^{l+k}\tbinom{k}{l}\Delta_{ul,j}G(x).
\]
\begin{definition}
\label{def:nikolskii}
For given  vectors $\vec{\beta}=(\beta_1,\ldots,\beta_d)\in (0,\infty)^d$, $\vec{r}=(r_1,$ $\ldots,r_d)\in [1,\infty]^d$,
 and $\vec{L}=(L_1,\ldots, L_d)\in (0,\infty)^d$ we
say that function $G:\bR^d\to \bR$ belongs to  anisotropic
Nikolskii's class $\bN_{\vec{r},d}\big(\vec{\beta},\vec{L}\big)$ if
 $\|G\|_{r_j}\leq L_{j}$ for all $j=1,\ldots,d$ and
 there exist natural number  $k_j>\beta_j$ such that
\[
 \big\|\Delta_{u,j}^{k_j} G\big\|_{r_j} \leq L_j |u|^{\beta_j},\;\;\;\;
\forall u\in \bR,\;\;\;\forall j=1,\ldots, d.
\]
\end{definition}
In addition to constraint
$f\in \bN_{\vec{r},d}(\vec{\beta},\vec{L})$ we also assume that
$f\in \bB_q(Q)$.
By  definition of Nikolskii's class,  $f\in\bN_{\vec{r},d}(\vec{\beta},\vec{L})$ implies
$f\in\bB_{r^*}(\max_{l=1,\ldots,d}L_l)$, where $r^*:=\max_{l=1,\ldots,d}r_l$.
Since we are interested in estimating  $\|f\|_p$,
it is necessary to suppose that this norm is bounded.
Therefore in all what follows we  assume that
$q\geq p\,\vee\, r^*$.
\par
Asymptotic behavior of the minimax risks on  anisotropic Nikolskii's classes is conveniently expressed in terms
of
the following parameters:
\begin{align*}
&
\frac{1}{\beta}:=\sum_{j=1}^d\frac{1}{\beta_j},\quad \frac{1}{\omega}:=\sum_{j=1}^d\frac{1}{\beta_jr_j},\quad
\BL:=\prod_{j=1}^dL_j^{1/\beta_j},
\\
&
\tau(s):=1-\frac{1}{\omega}+\frac{1}{\beta s},\;\;\; s\in[1,\infty].
\end{align*}
It is worth mentioning that quantities $\tau(\cdot)$ appear in embedding theorems for Nikolskii's classes; for details
see
\cite{Nikolskii}.
\par
Now we
are ready to state lower bounds on the minimax risk in the problem of estimating
the $\bL_p$--norm $\|f\|_p$. We consider the cases of integer and non--integer $p$ separately.
\subsection{The case of integer $p\geq 2$}
Define
\begin{gather*}
\theta:=\left\{
\begin{array}{clc}
\frac{1}{\tau(1)},\quad&\tau(p)\geq 1;
\\*[2mm]
\frac{1/p-1/q}{1-1/q-(1-1/p)\tau(q)},\quad&\tau(p)< 1,\;\tau(q)<0;
\\*[2mm]
\frac{\tau(p)}{\tau(1)},\quad&\tau(p)<1,\;\tau(q)\geq 0,
\end{array}
\right.
\end{gather*}
and let
\[
\phi_n:=\BL^{\frac{1-1/p}{\tau(1)}}n^{-\theta^*},\quad \theta^*:=2^{-1}\wedge \theta.
\]
\begin{theorem}
\label{th:lb-integer-r}
For any $\vec{\beta}\in (0,\infty)^d$, $\vec{L}\in(0,\infty)^d$, $\vec{r}\in[1,\infty]^d$, $q\geq p\vee r^*$ and
$p\in\bN^*$, $p\geq 2$
there exists  $c>0$ independent of $\vec{L}$ such that
\begin{equation*}
\liminf_{n\to\infty} \phi^{-1}_n\, \cR_n\big[\bN_{\vec{r},d}\big(\vec{\beta},\vec{L}\big)\cap \bB_q(Q)\big]\geq c.
\end{equation*}
\end{theorem}
\begin{remark}\label{rem:consistency}
In the companion paper \cite{GL20b}
we demonstrate that the rates of convergence of the minimax risk established
in Theorem~\ref{th:lb-integer-r}  are  {\em minimax}, that is,  they are attained by explicitly
constructed estimation procedures.
\end{remark}
\par
The lower bounds on the minimax rates of convergence of Theorem~\ref{th:lb-integer-r}
exhibit rather unusual features as compared to the results on estimating the $\bL_p$--norm of a signal
in the Gaussian white noise model [see \cite{LNS} and \cite{han&co2}].
\par\smallskip
1. It is quite surprising that the obtained asymptotics of the minimax risk does not depend on $p$ and $q$ if
$\tau(p)\geq 1$. Perhaps it is even more surprising that
in some cases the $\bL_p$-norm of a probability density can be estimated with the {\em parametric rate}!
On the other hand, it is easily seen that $\theta <1/2$ if $\tau(p)< 1,\;\tau(q)<0$; therefore, the parametric
rate is not achievable in this regime.
\par \smallskip
2.  If $r^*=\max_{l=1,\ldots, d} r_l \leq p$ and $q=p$ then
uniformly consistent estimators over anisotropic Nikol'skii's classes do not exist when $\tau(p)\leq 0$. This
together with Remark~\ref{rem:consistency} implies that condition $\tau(p)>0$ is
{\em necessary and sufficient} for  existence of  uniformly consistent estimators of the $\bL_p$-norm.
\par\smallskip
3. Taking together the previous remarks, we see that
in the considered estimation problem
the full spectrum of asymptotic behavior for the minimax risk
is possible: from parametric rate of convergence to inconsistency.
To the best of our knowledge this phenomenon has not  been  observed before.
\subsection{The case of  non-integer $p\geq 1$}
Define
\begin{eqnarray*}
\vartheta &:=&\left\{
\begin{array}{clc}
\frac{1}{2}\wedge \frac{1-1/p}{\tau(1)},\quad&\tau(p)\geq 1-1/p;
\\*[2mm]
\frac{1/p-1/q}{1-1/q-\tau(q)},\quad&\tau(p)< 1-1/p,\;\tau(q)<0;
\\*[2mm]
\frac{1}{2}\wedge\frac{\tau(p)}{\tau(1)},\quad&\tau(p)<1-1/p,\;\tau(q)\geq 0,
\end{array}
\right.
\end{eqnarray*}
\begin{eqnarray}\label{eq:vartheta*}
\vartheta^*&:=&
\left\{
\begin{array}{cl}
\frac{2(1-1/p)-\tau(p)}{\tau(1)},\; & \; \tau(p)\geq 1-1/p;
\\*[2mm]
\vartheta,\; &\;\tau(p)< 1-1/p,\;\tau(q)<0;
\\*[2mm]
2p,\; &\;\tau(p)<1-1/p,\;\tau(q)\geq 0,
\end{array}
\right.
\end{eqnarray}
and let
$$
\phi_n:=\BL^{\frac{1-1/p}{\tau(1)}}n^{-\vartheta}\big[\ln(n)\big]^{\vartheta^*-2p}.
$$
\begin{theorem}
\label{th:lb-noninteger-r}
For any $\vec{\beta}\in (0,\infty)^d$, $\vec{L}\in (0, \infty)^d$,  $\vec{r}\in[1, \infty]^d$
and $p\notin\bN^*, p>1$ there exists $c>0$ independent of $\vec{L}$ such that
\begin{equation*}
\liminf_{n\to\infty}\,\phi^{-1}_n\,\cR_n\big[\bN_{\vec{r},d}\big(\vec{\beta},\vec{L}\big)\cap \bB_q(Q)\big]\geq c.
\end{equation*}
\end{theorem}
\par
1.~Note that
the rates  of convergence established in Theorems \ref{th:lb-integer-r} and \ref{th:lb-noninteger-r} are different, except the case $\tau(p)<1-1/p,\;\tau(q)\geq 0$.
As we can see, the estimation accuracy for integer values of  $p$
is much better than for the non--integer ones. For the first time this phenomenon was discovered by
 \cite{LNS} in the problem of estimating the $\bL_p$--norm of a signal
 in  the univariate Gaussian white noise model.
\par
2. 
Theorem~\ref{th:lb-noninteger-r} shows
that if $q=p$ and $\tau(p)=0$ then there
no uniformly consistent
estimator exists. If $q=p$ and $\tau(p)<0$ then the lower bound becomes logarithmic in~$n$. We
conjecture that in this case  there are no uniformly consistent
estimators as well. If our conjecture is true, the proof of lower bounds will require
some additional considerations.
\par
3. It is not difficult to check that
the rate of convergence corresponding to the zone $\tau(p)\leq 1-1/p$ is slower than the one corresponding
to $\tau(p)>1-1/p$ independently of the value of $q$.
\subsection{How optimal are the  risk bounds of Theorem~\ref{th:lb-noninteger-r}?}
As it was mentioned above,
in this paper we do not discuss estimation procedures; we refer to  \cite{GL20b}
for construction of rate--optimal estimators of~$\|f\|_p$ for integer values of~$p$.
However, for non--integer values of~$p$ in
some cases very simple constructions
lead to {\em nearly rate--optimal adaptive estimators} of $\bL_p$-norms.
Let us discuss one such estimator under condition that density $f$ is uniformly bounded, i.e.,~$q=\infty$.
\par
Let $p\notin\bN^*, p>1$  be fixed.
Remark that if $q=\infty$ then
\begin{align*}
&\vartheta=\left\{
\begin{array}{clc}
\frac{1}{2}\wedge \frac{1-1/p}{\tau(1)},\quad&\tau(p)\geq 1-1/p;
\\*[2mm]
\frac{\omega}{p},\quad&\tau(p)< 1-1/p,\;\tau(\infty)<0,
\end{array}\right.
\\*[2mm]
&\phi_n=\BL^{\frac{1-1/p}{\tau(1)}}n^{-\vartheta}\big[\ln(n)\big]^{\vartheta-2p}.
\end{align*}
Consider the following sets of parameters:
\begin{eqnarray*}
\mD_1&=&\big\{\big(\vec{\beta},\vec{r}\big):\; \tau(p)> 2(1-1/p)\big\};
\\
\mD_2&=&\big\{\big(\vec{\beta},\vec{r}\big):\; \tau(p)< 1-2/p,\; \tau(\infty)<0\big\}.
\end{eqnarray*}
Let $\ell>0$ be an arbitrary a~priori chosen real  number, and  let
$$
\varphi_n:=\left\{
\begin{array}{cc}
(\BL/n)^{\frac{1-1/p}{\tau(1)}}\big[\ln(n)\big]^{d-1+\frac{1-1/p}{\tau(1)}},\; & (\vec{\beta},\vec{r})\in\mD_1;
\\*[2mm]
(\BL\ln(n)/n)^{\omega/p},\; & (\vec{\beta},\vec{r})\in\mD_2.
\end{array}
\right.
$$
\par
Let
$\widehat{f}(x)$, $x\in\bR^d$ be the estimator of $f(x)$
built in Theorem~1 of \cite{LW} in the case $\alpha=0$ [see also \cite{gl14}], and consider {\em
the plug--in estimator}
of the $\bL_p$--norm, $\widehat{F}:=\|\widehat{f}\|_p$.
\begin{theorem}
\label{th:ub-noninteger-r}
For any $Q>0$, $ L_0>0,\ell >0$, $\vec{L}\in [L_0,\infty)^d$ and any $\vec{\beta}\in (0,\ell]^d$, $\vec{r}\in (1,\infty]^d$ belonging to $\mD_1\cup\mD_2$
there exists $C<\infty$, independent of $\vec{L}$,  such that
\begin{equation*}
\limsup_{n\to\infty}\,\varphi^{-1}_n\,\cR_n\big[\widehat{F};\; \bN_{\vec{r},d}\big(\vec{\beta},\vec{L}\big)\cap\bB_\infty(Q)\big]\leq C.
\end{equation*}
\end{theorem}
The proof of this theorem is trivial. By the triangle inequality $|\widehat{F}-\|f\|_p|\leq \|\widehat{f}-f\|_p$, so
that the problem of estimating $\|f\|_p$ can be reduced to the problem of adaptive
estimation of $f$
under the $\bL_p$-loss. The stated upper bound follows from the results of Theorem~3
in \cite{LW} corresponding
to what is called in that paper  {\em tail zone} and {\em sparse zone~1}.
Combining bounds of Theorems~\ref{th:lb-noninteger-r} and~\ref{th:ub-noninteger-r}
we come to the following statement.
\begin{corollary}
\label{assert:asympt-minimax-noninteger-r}
For any $Q>0$, $ L_0>0,\ell >0$, $\vec{L}\in [L_0,\infty)^d$ and any $\vec{\beta}\in (0,\ell]^d$, $\vec{r}\in (1,\infty]^d$ belonging to $\mD_1\cup\mD_2$
one has for all $n$ large enough
\begin{equation*}
\big[\ln(n)\big]^{\gamma_1 - 2p}\;\lesssim\; n^{\vartheta}\,
\cR_n\big[\bN_{\vec{r},d}\big(\vec{\beta},\vec{L}\big)\cap \bB_\infty(Q)\big] \;\lesssim\; \big[\ln(n)\big]^{\gamma_2},
\end{equation*}
where
\[
 \gamma_1:=\left\{\begin{array}{cl}
                   \frac{2(1-1/p)-\tau(p)}{\tau(1)}, & \mD_1,\\
                   \omega/p, & \mD_2,
                  \end{array}
\right.\;\;\;\;
\gamma_2:=\left\{\begin{array}{cl}
                   d-1+\frac{1-1/p}{\tau(1)}, & \mD_1,\\
                   \omega/p, & \mD_2.
                  \end{array}
\right.
\]
\end{corollary}
Thus the estimator  $\widehat{F}$ is {\em nearly rate--optimal adaptive} over the scale of Nikolskii's classes whose parameters belong to  $\mD_1\cup\mD_2$.
\subsection{Relation to  density estimation under $\bL_p$-loss}
Assume that we are interested in estimating density $f$ under $\bL_p$-loss from the
observation $X^{(n)}$. Measuring accuracy of estimation procedures by the $\bL_p$-loss leads
to the quadratic risk in the form
$$
\mR_n[\widetilde{f}, f]:=\Big(\bE_f \big[\|\widetilde{f}-f\|_p\big]^2\Big)^{1/2}.
$$
In view of the triangle inequality we have for any $f\in\mF$
$$
\mR_n[\widetilde{f}, f]\geq \cR_n[\widetilde{F}, f],
$$
where $\widetilde{F}=\|\widetilde{f}\|_p$. Hence, whatever the functional class $\bF$ is,
one has
$$
\mR_n[\bF]:=\inf_{\widetilde{f}}\sup_{f\in\bF}\mR_n[\widetilde{f}, f]\geq \cR_n[\bF].
$$
and we assert that any lower bound for $ \cR_n[\bF]$ is automatically the lower bound for $\mR_n[\bF]$.
In particular, assuming that $r^*\leq p$ and putting $q=p$ we deduce from Theorems \ref{th:lb-integer-r}
and \ref{th:lb-noninteger-r}
the following result.
\begin{corollary}
\label{cor:L_p-loss}
Let either $0<\tau(p)< 1-1/p$ if $p\notin\bN^*$ or $0<\tau(p)< 1$ if $p\in\bN^*$. Then
\begin{equation*}
\liminf_{n\to\infty}\,n^{-\frac{\tau(p)}{\tau(1)}}\,\mR_n\big[\bN_{\vec{r},d}\big(\vec{\beta},\vec{L}\big)\big]\geq c>0.
\end{equation*}
\end{corollary}
This  asymptotics of the minimax risk
is related to the fact that the estimated density may be unbounded.
To the best of our knowledge, this result is new.
In the one-dimensional case for $p=2$ the obtained rate    coincides
with the one in \cite{Birge}.

\section{Lower bounds for estimation of  general non--linear functionals}
\label{sec:abstract-LBT}
The results of Theorems~\ref{th:lb-integer-r} and~\ref{th:lb-noninteger-r} follow from
general machinery
for derivation
of lower bounds on minimax risks in the density model.
In this section we
develop  this technique in full generality for a broad class of nonlinear functionals to be estimated.
\par
Let $G:\bR\to\bR$ and $H:\bR_+\to\bR$ be fixed functions. We are interested in estimating
the functional
\begin{equation}\label{eq:Psi-G-H}
\Psi(f)=G\bigg(\int_{\bR^d}H\big(f(x)\big)\rd x\bigg)
\end{equation}
from observation $X^{(n)}=(X_1,\ldots, X_n)$.
Let $\bF$ be a  class of  functions defined on $\bR^d$ and  let
$$
\cR_n\big[\bF\big]:=\inf_{\tilde{\Psi}}\sup_{f\in\bF}\Big(\bE_f \big[\tilde{\Psi}-\Psi(f) \big]^2\Big)^{1/2},
$$
where the infimum is taken over all possible estimators of $\Psi$.
Our goal is to derive an explicit lower bound on the minimax risk
under mild condition on functions $G$ and $H$ and functional class
$\bF$.
\par
We remark that the class of considered functionals is rather broad and includes many problem instances of interest.
Let us give some examples.
\begin{enumerate}
\item
Let  $G(y)=y^{1/p}$ and $H(y)=y^p$ for some $p\in(1,\infty)$; then $\Psi(f)$ is  the $\bL_p$-norm of~$f$, and
estimation of this functional is the subject of the present paper.
\item
The choice $G(y)=ay$ leads to the estimation of the {\em integral-type functionals}.
The following particular cases have been considered in the literature.
\begin{itemize}
\item[(a)] if  $a=1$ and $H(y)=y^p$ with $p\in\bN^*, p\geq 2$, then the corresponding  functional is
$\Psi(f)=\|f\|_p^p$; see for instance, \cite{bickel-ritov}, \cite{picard},  \cite{beatrice}, \cite{waart};
\item[(b)]
the case $a=1$ and $H(y)=-y\ln(y)$ corresponds to {\em the differential  entropy},  $\Psi(f)=-\int f(x)\ln f(x)\rd x$; see, e.g., \cite{koza};
\item[(c)]
if  $a=(p-1)^{-1}$ and $H(y)=y-y^p$ with $ p\neq 1$ then $\Psi(f)$ is {\em the Tsallis  entropy},
$\Psi(f)=(p-1)^{-1}(1-\int |f(x)|^p\rd x)$;
see \cite{tsallis}, \cite{leon}.
\end{itemize}
\item Let $G(y)=(1-p)^{-1}\ln(y)$ and $H(y)=y^p$ with $p\neq 1$; then the corresponding functional is
{\em the R\'enyi entropy}, $\Psi(f)=(1-p)^{-1}\ln\big(\int |f(x)|^p\rd x\big)$; see \cite{renyi}, \cite{leon}.
\end{enumerate}
\par
The technique for derivation of lower bounds relies upon construction of a parameterized
family of functions equipped with a pair prior probability measures on it. Below
we discuss  these
construction ingredients in succession.

\subsection{Parameterized family of functions}
Let $\Lambda:\bR^d\to \bR_+$  be a function satisfying the following conditions:
\begin{equation}
\label{eq:Lambda}
\Lambda(x)=0,\; \; \forall x\notin [-1,1]^d,\quad \int_{\bR^d}\Lambda(x)\rd x =1.
\end{equation}
 Let $|\cdot|_\infty$ denote the $\ell_\infty$--norm on $\bR^d$, and let $\cM$ be a given finite set of indices  of cardinality  $M={\rm card}(\cM)$.
Let
 $\big\{x_{m}\in\bR^d,\;m\in \cM\big\}$ be a finite set of points in $\bR^d$
 satisfying
$$
 \big|x_{k}-x_{m}\big|_\infty\geq 2,\quad\;\forall k\neq m,\; k, m\in\cM.
$$
Fix vector $\vec{\sigma}=(\sigma_1,\ldots,\sigma_d)\in (0,1]^d$ and constant $A>0$
and define for any $m\in\cM$
\begin{equation}\label{eq:Pi-m}
\Lambda_m(x)=A\Lambda\big((x-x_{m})/\vec{\sigma}\big),\quad \Pi_{m}=
\Big\{x\in\bR^d:\; \big|(x-x_{m})/\vec{\sigma}\big|_\infty\leq 1\Big\},
\end{equation}
where the division is understood in the coordinate--wise sense.
In words, $\Pi_m$ is a rectangle in $\bR^d$ centered at $x_m$ with edges of half--lengths $\sigma_1, \ldots \sigma_d$
that are parallel to the coordinate axes.
It is obvious that $\Lambda_m$ is supported on $\Pi_m$ for any  $m\in\cM$, and $\Pi_m$ are disjoint:
\begin{eqnarray}
\label{eq5:proof-th:lower-bound-deconvolution}
&& \Pi_m\cap \Pi_k=\emptyset,\quad \forall k\neq m,\; k,m\in\cM.
\end{eqnarray}
Let $\Pi_{0}:=\bR^d\setminus\cup_{m\in\cM}\Pi_{m}$, $\bma:=\prod_{l=1}^d\sigma_l$, and
\begin{eqnarray*}
&&
\varrho_w(z) :=\sum_{m\in\cM}\big(w_{m}\big)^z,\quad w\in[0,1]^{M},\;\;z>0.
\end{eqnarray*}
\par
Let $f_{0}$ be a probability density supported on $\Pi_{0}$.
Define  the family of functions:
 \begin{equation}\label{eq:f-w}
 f_{w}(x):=\big[1-A\bma\varrho_{w}(1)\big]f_{0}(x)+A\sum_{m\in\cM}w_{m}\Lambda_m(x), \quad
 w\in[0,1]^{M}.
 \end{equation}
The family $\{f_w, w\in [0,1]^M\}$ involves tuning
parameters $A$, $\vec{\sigma}$ and $M$ that will be specified in the sequel.
 The most important element of our approach consists in equipping  $[0,1]^M$ with
 two product probability measures, thus
 assuming that~$w$ is a random vector distributed in accordance with one of them.
 Then functions $f_w$ become random, and they are not necessarily
 density functions and/or functions from the functional class $\bF$ for all realizations of $w$.
 %
 With  conditions introduced below we ensure that $f_w\in \bF\cap \mF$ with
 large enough probability.

\subsection{Prior probability measures}
Let
$\mP[0,1]$ be  the set of all probability measures with total mass on $[0,1]$.
For  $\mc \in\mP[0,1]$,  $z\geq 0$  we define
\begin{equation*}
e_{\mc}(z):=\int_{0}^1 x^z \mc(\rd x).
\end{equation*}
\par
\begin{definition}
\label{def1:prior}
For a pair of probability measures $\ma, \mb\in\mP[0,1]$ we will write
$\ma\stackrel{}{\backsim}\mb$ if    $e_\ma(1)=e_\mb(1)$.
\end{definition}
\par
Let  $\bze:=(\zeta_{m}, m\in\cM)$ be independent identically distributed random variables, and $\zeta_m$ is distributed
$\mc\in \mP[0,1]$, $m\in \cM$.
 The law of $\bze$ and the corresponding expectation
 will be denoted by $\bP_{\mc}$ and $\bE_{\mc}$ respectively.
Define
\[
 \mmp_{\bze}(x) := \prod_{i=1}^d f_{\bze}(x_i),\;\;\;x\in \bR^{dn};
\]
here and from now on we regard $x=(x_1, \ldots, x_n)$, $x_i\in \bR^d$ as an element
of $\bR^{dn}$.
\subsection{Assumptions on the family of functions and prior measures}
Now we introduce general assumptions that relate properties of
parameterized family
$\{f_w, w\in [0,1]^M\}$ and prior measures on $[0,1]^M$.
\begin{assumption}
\label{ass1:mixed}
There exist $\e\in (0,1)$ and  two probability measures $\ma,\mb\in\mP[0,1]$, $\ma\backsim \mb$
such that
$$
\bP_{\mc}\big\{f_{\bze}\in \bF\big\}\geq 1-\e, \quad \mc\in\{\ma,\mb\},
$$
where $f_{\bze}$ is defined in (\ref{eq:f-w}).
\end{assumption}
\par
Assumption \ref{ass1:mixed} stipulates that
under  prior probability measures $\ma$ and $\mb$ random function $f_{\bze}$
belongs to the functional class $\bF$ with probability at least $1-\e$.
Note that this assumption does not guarantee that $f_{\bze}$  is a probability density;
by construction, only assumption $\int f_{\bze}=1$ is fulfilled   for all realizations of $\bze$.
\par
We also need conditions
that relate parameters $A,\vec{\sigma}, M$
of the family of functions
with the number of observations $n$.
\begin{assumption}
\label{ass2:mixed}
For sufficiently small  $\kappa_0>0$ and  sufficiently large $\upsilon>0$  one has
\begin{eqnarray}
\label{eq5-lbf}
 A\bma\sqrt{M}&\leq& \kappa_0n^{-1/2},
\\*[2mm]
\label{eq3**-lbf}
 M&\geq& 36\upsilon,
\\*[2mm]
A\bma M&\leq&  1/2.
\label{eq3-lbf}
\end{eqnarray}
\end{assumption}
Condition (\ref{eq3**-lbf})
guarantees that random function $f_{\bze}$ concentrates properly around its expectation, while
(\ref{eq3-lbf}) implies that $f_{\bze}$ is a probability density for all realizations of $\bze$.
Indeed, by construction $\int f_{\bze}(x) \rd x=1$, and, in view of (\ref{eq3-lbf}),  $f_{\bze}\geq 0$
because  $\varrho_{\bze}(1)\leq M$ for all~$\bze$.
In addition,  condition (\ref{eq5-lbf}) allows us to construct  a product form
approximation for the Bayesian likelihood ratio
$\bE_{\ma} [\mmp_{\bze}(\cdot)]/\bE_{\mb}[\mmp_{\bze}(\cdot)]$, which is an essential step in
the  derivation of  lower bounds.
\subsection{Main results}
To state  lower bounds on the minimax risk for estimating functional $\Psi(f)$
we require  notation that involves functions $H$ and $G$ appearing in (\ref{eq:Psi-G-H}).
\par
Define functions
$S_0:[0,1]\to \bR$ and $S:\bR_+\to \bR$ by
\begin{align}\label{eq:S0-S}
S_0(z) :=\int_{\Pi_{0}}H\big((1-z)f_{0}(x)\big)\rd x,  \;\;\;\;
S(z):=\int_{[-1,1]^d}H\big(z\Lambda(x)\big)\rd x.
\end{align}
For $\mc\in\mP[0,1]$ let
\begin{align}\label{eq:E-pi-V-pi}
 E_{\mc}(A) :=\int_{0}^1S(Ay)\mc(\rd y), \;\;\;\;
 V_{\mc}(A) :=\bigg[\int_{0}^1S^2(Ay)\mc(\rd y)\bigg]^{1/2}.
\end{align}
We tacitly assume that $E_{\mc}(A)$ and $V_{\mc}(A)$ are finite for all $A>0$ and for all considered
measures $\pi\in \mP[0,1]$; fulfillment of  this assumption should be verified for
every concrete problem instance.
\par
To clarify the meaning of notation
introduced in (\ref{eq:S0-S}) and (\ref{eq:E-pi-V-pi}) we observe that for  the given family of functions
$\{f_w, w\in [0,1]^M\}$ one has
\begin{align}\label{eq:int-H}
 \int_{\bR^d} H\big(f_{\bze}(x)\big)\rd x
 &= S_0(A\bma \varrho_{\bze}(1)) + \bma \sum_{m\in \cM} S(A\zeta_m)
 \nonumber
 \\
 &\approx
 S_0(A\bma Me_\mc(1))+ \bma M E_\mc (A),
\end{align}
where the approximate equality in the second line designates that
the sums of independent random variables $\varrho_{\bze}(1)$ and $\sum_{m\in \cM} S(A\zeta_m)$
concentrate properly around their expectations $Me_\mc(1)$ and
$ME_{\mc}(A)$ respectively.
In addition, the lower bound  derivation requires analysis of
discrepancy between the values
of the functional $\Psi(f_{\bze})$ when $\bze$ is distributed according to prior measures
$\ma, \mb\in \mP[0,1]$, $\ma\backsim \mb$.
This fact along with (\ref{eq:int-H})   motivates the following notation.
\par
Let $\ma, \mb \in \mP[0,1]$, $\ma\backsim \mb$ and $\mc\in \{\ma, \mb\}$. Define
\begin{eqnarray}
 H^*_{\mc}&:=& S_0(A\bma M e_\mc(1))+ \bma ME_\mc(A),\;\;
 \label{eq:H-pi}
 \\*[2mm]
 \alpha_{\mc}&:=&\eta_{S_0}\big(A\bma M e_\mc(1);\, A\bma \sqrt{M}\big)+  \bma \sqrt{\upsilon M}
 \,V_{\mc}(A),
 \label{eq:alpha-pi}
 \\*[2mm]
 \cJ_{\mc} &:=& \Big[\inf_{|\alpha|\leq \alpha_{\mc}} G(H^*_{\mc}+\alpha),
 \sup_{|\alpha|\leq \alpha_{\mc}} G(H^*_{\mc}+\alpha)\Big],
\label{eq:cJ-pi}
 \end{eqnarray}
where $\eta_{S_0}(x;\delta)$ stands for the local modulus of continuity of function $S_0$,
\[
 \eta_{S_0}(x; \delta):= \sup_{y:|y-x|\leq \delta} |S_0(x)-S_0(y)|,\;\;x, y\in [0,1],\;\;\delta>0.
\]
Define also
\[
 \Delta(\ma, \mb):= \min\big\{|x-x^\prime|: x\in \cJ_{\ma}, \, x^\prime\in \cJ_{\mb}\big\};
\]
clearly, $\Delta(\ma, \mb)$ is the Hausdorff  distance between the intervals
$\cJ_{\ma}$ and $\cJ_{\mb}$.
\par
Finally we let
\begin{equation*}
n_{m}(x):=\sum_{i=1}^n\mathrm{\rm 1}_{\Pi_{m}}(x_i),\;\;
n_{0}(x):=n-\sum_{m\in\cM}n_{m}(x),\
\;\;x\in\bR^{dn},
\end{equation*}
where sets $\Pi_m$, $m\in \cM$ are defined in  (\ref{eq:Pi-m}).
The quantities $n_m(x)$ and $n_0(x)$ have evident probabilistic interpretation: if
$X^{(n)}=(X_1, \ldots, X_n)$ is the sample then
$n_m(X^{(n)})$ and $n_0(X^{(n)})$ are numbers of observations
in the sets $\Pi_m$ and $\Pi_0$ respectively.
Furthermore, for a pair of measures $\ma\backsim \mb$ we define
\begin{align}
&\Upsilon(x):=
\prod_{m\in\cM}\frac{\gamma_{m,\ma}(x)}{\gamma_{m,\mb}(x)}, \;\;\; x\in\bR^{dn},
\nonumber
\\
&\gamma_{m,\mc}(x):= \int_{0}^1 y^{n_{m}(x)}e^{-Dn_{0}(x)y}\mc(\rd y),\;\;\;D:= \frac{A \bma}{1- A \bma Me_\mc(1)}.
\label{eq:bayes-ration-second2}
\end{align}
Observe that $D$
does not depend on $\pi\in \{\ma, \mb\}$ because $\ma\backsim \mb$.
\par
Now we are in a position to formulate the main result of this section.
\begin{theorem}
\label{th:lower-bound-general}
Let $\ma, \mb\in \mP[0,1]$, $\ma\backsim \mb$, and suppose that
Assumptions \ref{ass1:mixed} and \ref{ass2:mixed}  are  fulfilled.
If $\Delta(\ma, \mb)>0$ then
\begin{align}
 \;\;\;\;36e [\Delta (\ma, & \mb)]^{-2}\,  \cR^2_n[\bF]
 \nonumber
 \\*[2mm]
  &\geq
\bE_{\mb}\Big[
\bP_{f_{\bze}}\big\{\Upsilon\big(X^{(n)}\big)\geq \tfrac{1}{2}\big\}\Big]
-\upsilon^{-1}
-  4\sqrt{2(2\upsilon^{-1}+\e)}.
\label{eq:lower-bound}
\end{align}
\end{theorem}
In order to apply general
lower bound (\ref{eq:lower-bound}) in concrete problem instances we need  to compute or bound from below
the quantity $\Delta(\ma, \mb)$ and to show that  the right hand side is strictly positive.
The next
two corollaries of Theorem~\ref{th:lower-bound-general}
derive lower bounds on the right hand side of (\ref{eq:lower-bound})
under additional conditions on prior measures and parameters of the family $\{f_w, w\in [0,1]^M\}$.
\begin{definition}
\label{def2:prior}
Let $r\in\bN^*, r\geq 1$ be fixed.
For $\ma,\mb\in\mP[0,1]$ we write  $\ma\stackrel{r}{\backsim}\mb$
if $e_{\ma}(k)=e_{\mb}(k)$ for all  $k=1,\ldots, r$.
\end{definition}
Proposition~\ref{prop2} of Section~\ref{sec:prior-construction} presents a
sophisticated  construction of probability measures $\ma$, $\mb$ satisfying
conditions of Definition~\ref{def2:prior} and possessing some additional properties.
\begin{corollary}
\label{cor:lower-bound-genera2}
Let $r$ be a positive integer number, possibly dependent on $n$, such that
$r>\ln(36\upsilon)$,
and let $\ma, \mb\in \mP[0,1]$ satisfy $\ma\stackrel{2r}{\backsim} \mb$. Suppose that
Assumptions \ref{ass1:mixed} and \ref{ass2:mixed}  are  fulfilled.
Assume that for sufficiently small $\kappa_1>0$ one has
\begin{eqnarray}
\label{eq:condition-main2}
&& n\mA\bma \leq \kappa_1r,\;\;\;\; M\leq e^{r}.
\end{eqnarray}
If  $\Delta(\ma,\mb)>0$ then putting $C_*=(36e)^{-1/2}\big[3^{-1} - \upsilon^{-1} -\sqrt{2(\e+2/\upsilon)}\,\big]^{1/2}$ one has
$$
\cR_n[\bF]\geq  C_*\Delta(\ma,\mb).
$$

\end{corollary}
\begin{definition}
\label{def3:prior}
Let $t,r\in\bN^*$, $r\geq t>1$ be fixed.
For $\ma,\mb\in\mP[0,1]$ we write that  $\ma\stackrel{r,t}{\backsim}\mb$
if the following requirements are fulfilled:
\begin{enumerate}
\item $e_{\ma}(k)=e_{\mb}(k)$ for all  $k=1,\ldots, r,\; k\neq t$;
\smallskip
\item $e_{\ma}(t)\neq e_{\mb}(t)$.
\end{enumerate}
\end{definition}
Proposition~\ref{prop1} of Section~\ref{sec:prior-construction} presents a construction of measures $\ma$ and $\mb$
satisfying $\ma\stackrel{r,t}{\backsim}\mb$.
\begin{corollary}
\label{cor:lower-bound-genera3}
Fix $r\geq t>1$, and  let $\ma, \mb\in \mP[0,1]$ satisfy $\ma\stackrel{r, t}{\backsim} \mb$.
Suppose that Assumptions~\ref{ass1:mixed}--\ref{ass2:mixed} are  fulfilled, and
for sufficiently small $\kappa_1>0$ and $t$ independent of $n$ one has
\begin{eqnarray}
\label{eq:condition-main}
 n A\bma M^{1/t} \,\leq\, \kappa_1.
\end{eqnarray}
If   $\Delta(\ma,\mb)>0$   then
$$
\cR_n[\bF]\geq C_*\Delta(\ma,\mb),
$$
where $C_*$ is defined in Corollary~\ref{cor:lower-bound-genera2}.
\end{corollary}
\subsection{Discussion}
In this section we discuss  applicability of Theorem~\ref{th:lower-bound-general} and Corollaries~\ref{cor:lower-bound-genera2} and~\ref{cor:lower-bound-genera3}, and   main ideas that underlie the proofs of these results.
\paragraph*{More general statistical experiments}
The proofs of Theorem \ref{th:lower-bound-general} and
Corollaries~\ref{cor:lower-bound-genera2} and~\ref{cor:lower-bound-genera3} do not use
the fact that density $f$ is defined on~$\bR^d$. In fact, after minor changes and
modifications our construction  is applicable  in an arbitrary density model.
\par
Let $(\cX,\mB,\lambda)$ be a measurable space, and let $X$ be an
$\cX$-valued random variable whose law has the density $f$ with respect to   measure $\lambda$.
Assume that we observe $X^{(n)}=(X_1,\ldots, X_n)$, where $X_i, i=1,\ldots,n$, are independent copies of $X$. The goal is to estimate the functional
\begin{equation*}
\Psi(f)=G\bigg(\int_{\cX}H\big(f(x)\big)\lambda(\rd x)\bigg),
\end{equation*}
where as before  $G:\bR\to\bR$ and $H:\bR_+\to\bR$  are fixed functions.
\par
Let $\cM$ be a finite set of indices with cardinality $M$, possibly dependent on $n$.
Let  $f_{0},\Lambda_m:\bX\to\bR_+$ and $\Pi_m\in\mB$, $m\in\cM$ be
collections of measurable functions  and sets satisfying the following conditions.
\begin{enumerate}
\item[(a)] $\Pi_m\cap\Pi_k=\emptyset$ for any $m\neq k, m,k\in\cM$;

\vskip0.2cm

\item[(b)]  $\lambda(\Pi_m)=\bma>0$ for any $m\in\cM$;

\vskip0.2cm

\item[(c)]  $\Lambda_m(x)=0, x\notin \Pi_m$ for any $m\in\cM$;

\vskip0.2cm

\item[(d)]  $\int_{\Pi_m}\Lambda_m(x)\lambda(\rd x)=1$ for any $m\in\cM$;

\vskip0.2cm

\item[(e)]  $\int_{\cX}f_0(x)\lambda(\rd x)=1$ and $f_0(x)=0$ for any $x\in\cup_{m\in\cM}\Pi_m$.

\end{enumerate}
Under these conditions
some evident minor modifications in definitions
should be made; for instance, function $S$ should be defined as
$$
S(z)=M^{-1}\sum_{m\in\cM}\int_{\cX}H\big(z\Lambda_m(x)\big)\lambda(\rd x).
$$
With these changes  the results of  Theorem~\ref{th:lower-bound-general}
and its corollaries remain valid.
\paragraph*{Method of proof}
The following
fundamental principles and main ideas  lie at the core
of
the proof
of Theorem~\ref{th:lower-bound-general} and its corollaries.
\par
The first idea
goes back to the paper \cite{LNS}. It reduces the original estimation problem to a problem of
testing two composite hypotheses for mixture distributions which are obtained
by imposing prior
probability measures with intersecting supports on  parameters of a functional family.
In~\cite{Tsybakov}
this technique is called  {\em the method of two fuzzy hypotheses}.
The choice of the prior measures  is based on {\em the moment matching technique};
see, e.g., \cite{cai-low-2}
and~\cite{han&co1}, where further references can be found.
%
%
We clarify the moment matching technique in Propositions~\ref{prop1} and~\ref{prop2} that can be
viewed as slight generalization and modification of the results in \cite{LNS}.
Detailed proofs of these statements are given in Appendix.
\par
The second idea is related to  construction of a specific parameterized family of densities on which the
lower bound on the minimax risks is established.
Here we use a construction that is similar to the one proposed in \cite{gl14}.
\par
The third main idea is related to the
analysis of the so-called {\em Bayesian likelihood ratio}.
This analysis, being common in estimating  nonlinear functionals,
depends heavily  on the considered statistical model.
The multivariate density model
on  $\bR^d$ requires development of an original technique
because,
in contrast to the Gaussian white noise or regression models,
the Bayesian likelihood ratio
is not
a product of independent random variables.
As a consequence,
standard methods based on computation of
the Kullback-Leibler, Hellinger or other divergences between distributions are not applicable.
That is why the proof of
Theorem~\ref{th:lower-bound-general} contains development
of two--sided  product--form bounds on the Bayesian likelihood ratio.
\subsection{Additional results}
In this section we discuss  implications of our results and the developed technique
for other problems of estimating nonlinear functionals.
In all examples below we consider functionals $\Psi(f)$ of type
(\ref{eq:Psi-G-H}) with $G(y)\equiv y$. Denote also
\[
\cF_{\Psi}:=\big\{f:\; |\Psi(f)|<\infty\big\}.
\]
\subsubsection*{Estimation of $\|f\|_p^p$, $p\in \bN^*$}
In this example $H(y)=y^p$, $p\in\bN^*$.
Apparently, the case  $p=2$ is  the most well studied setting.
Many authors, starting from the seminal  paper of \cite{bickel-ritov},
made fundamental contributions to the minimax and minimax adaptive estimation of quadratic functionals of probability density; see, for instance, \cite{birge-massart}, \cite{beatrice}
among many others.
\cite{picard} studied the case $p=3$  and
 \cite{waart} considered the setting with arbitrary integer $p$.
 It is worth noting that all aforementioned papers consider  either univariate or
 compactly supported densities belonging to a semi--isotropic functional class, that is $r_l=r$ for any $l=1,\ldots,d$.
\par
Let us consider the case $p=2$ and recall  one of the most well known  results.
Assume that the underlying density $f$ is {\em compactly supported} and belongs to the
anisotropic H\"older class $\bN_{\vec{\infty},d}\big(\vec{\beta},\vec{L}\big)$,
that is $r_l=\infty$  for all $l=1,\ldots,d$.
It is well known that tn this setting the minimax rate of convergence in estimating
$\|f\|_2^2$ is given by
\[
(1/n)^{\frac{4\beta}{4\beta+1}\wedge\frac{1}{2}};
\]
see, e.g., \cite{bickel-ritov} for one--dimensional case.
In particular, the parametric regime is possible if and only if
$\beta\geq 1/4$.
On the other hand, close inspection of  the proof of
Theorem \ref{th:lb-integer-r} shows that
the lower bound on the asymptotics of the minimax risk in estimating of $\|f\|_2^2$ is
 simply the squared rate
 found in this theorem   in the "nonparametric" regime.
 In particular, if $r_l=\infty$ for all $l=1,\ldots,d$ then Theorem~\ref{th:lb-integer-r}  yields
the rate
\[
(1/n)^{\frac{2\beta}{\beta+1}\wedge\frac{1}{2}}.
\]
We do not know whether this rate  is the minimax rate of convergence, but we can assert that
the parametric rate is not possible  if $\beta\leq 1/3$. This shows that problems of
estimating $\|f\|_2^2$  for compactly supported densities, and densities supported on the entire space
$\bR^d$ are completely different.
\par
Another interesting feature is that if $q=p$, $r^*\leq p$ and $\tau(p)\leq 0$ then
there is no uniformly consistent estimator of $\|f\|_2^2$ over anisotropic Nikolskii's class.
This phenomenon is again due
to the fact that the underlying density   is assumed to be supported  on the entire space $\bR^d$.
\subsubsection*{Estimation of the differential entropy}
This setting corresponds to
$H(y)=-y\ln(y)$. Applying the same reasoning as in the proof of Theorem~\ref{th:lb-noninteger-r}
in conjunction with Proposition~\ref{prop:timan} we are able  to prove the following  statement.
\begin{theorem}
\label{th:lb-entropy}
There exists $c>0$ such that for any $\vec{\beta}\in (0,\infty)^d$, $\vec{L}\in (0,\infty)^d$, $\vec{r}\in[1, \infty]^d$
\begin{equation*}
\liminf_{n\to\infty}\,[\ln n]^3\,\inf_{\widetilde{F}}\sup_{\cF_{\Psi}\cap\bN_{\vec{r},d}\big(\vec{\beta},\vec{L}\big)}\Big(\bE_f \big[\widetilde{F}-\Psi(f)\big]^2\Big)^{1/2}\geq c.
\end{equation*}
\end{theorem}
\subsubsection*{Estimation of $\|f\|_p^p$, $p\in (0,1)$} Let now  $H(y)=y^{p}, p\in (0,1)$.
The corresponding functional coincides up to a constant with
the Tsallis entropy with index  $p\in (0,1)$.
\begin{theorem}
\label{th:lb-Renyi}
For any $p\in (0,1)$ there  exists $c>0$ such that for any $\vec{\beta}\in (0,\infty)^d$.  $\vec{L}\in (0,\infty)^d$
and $\vec{r}\in[1, \infty]^d$
\begin{equation*}
\liminf_{n\to\infty}\, [\ln n]^{2p}\,\inf_{\widetilde{F}}\sup_{\cF_{\Psi}\cap\bN_{\vec{r},d}(\vec{\beta},\vec{L})}\Big(\bE_f \big[\widetilde{F}-\Psi(f)\big]^2\Big)^{1/2}\geq c.
\end{equation*}
\end{theorem}
The rates of convergence established in
Theorems~\ref{th:lb-entropy} and \ref{th:lb-Renyi} are very slow
and do not depend on the parameters
of the functional class. In particular, these results demonstrate    that smoothness  alone
is not sufficient in order to guarantee a  "reasonable" accuracy of estimation.
However, the recent paper \cite{han&co2} dealing with estimation of the differential  entropy shows
that
if underlying density satisfies moment conditions then the minimax risk converges to zero at the polynomial
(in $n$) rate. It is also clear that the   polynomial rate of convergence
in estimating of considered functionals is possible for smooth  compactly supported
densities.
The proof of Theorems \ref{th:lb-entropy} and \ref{th:lb-Renyi} coincides with the proof of Theorem~\ref{th:lb-noninteger-r}
up to minor modifications.
\section{Proofs of Theorems \ref{th:lb-integer-r} and \ref{th:lb-noninteger-r}}
\label{sec:proofs_Tms-norms}
We will prove only "nonparametric rates".
A lower bound corresponding to the parametric rate of convergence $n^{-1/2}$ can be easily derived
by reduction of the considered  problem to parameter estimation in a regular statistical model.
\par
\subsection{Preliminary remarks}
The proofs of Theorems \ref{th:lb-integer-r} and \ref{th:lb-noninteger-r} are based on application of Corollaries  \ref{cor:lower-bound-genera2} and~\ref{cor:lower-bound-genera3},
and contain many common elements.
In particular, in both proofs
we consider parameterized family of functions $\{f_w, w\in [0,1]^M\}$ defined in
(\ref{eq:Lambda})--(\ref{eq:f-w}), and we choose the sets $\Pi_0$ and $\Pi_m$, $m\in \cM$
so that
$\Pi_0\subset (-\infty, 0]^d$ and
$\Pi_m \subset [0, \infty)^d$ for all $m\in \cM$.
We equip the parameter set  $[0,1]^M$
with a pair of probability measures $\mu$, $\nu$
satisfying conditions of one of Definitions~\ref{def1:prior}--\ref{def3:prior}.
Along with  conditions of Definitions~\ref{def1:prior}--\ref{def3:prior} in the proofs of
Theorems~\ref{th:lb-integer-r} and~\ref{th:lb-noninteger-r}
we require that the specified probability measures $\mu$ and $\nu$ possess  the following property:
\begin{equation}
\label{eq:e>1/2}
 \sqrt{e_\pi (2z)}\leq 2e_\pi (z), \;\;\forall z\in\{r_1,\ldots,r_d,q, p\},\;\;\;\pi\in \{\mu, \nu\}.
\end{equation}
Recall that  $r_1,\ldots,r_d$ are the coordinates of the vector $\vec{r}$ used in the definition of the Nikolskii class.
Once measures satisfying conditions of Definitions~\ref{def1:prior}--\ref{def3:prior}
are constructed, they can be easily modified to satisfy (\ref{eq:e>1/2}); for details we
refer to the proofs of Propositions~\ref{prop1} and~\ref{prop2}. By convention, here and
from now on we put $[e_\pi(z)]^{1/z}=1$ for $z=\infty$.
\par
The
parameters
$A, \vec{\sigma}, M$ of the family  $\{f_w, w\in [0,1]^M\}$ are specified to guarantee
that under imposed prior measures
random functions $f_\bze$
 satisfy required
smoothness conditions with the probability controlled by parameter $\upsilon$. Under these circumstances,
parameter $\upsilon$
of Assumption~\ref{ass2:mixed} should be such that
constant $C_*$ in Corollaries~\ref{cor:lower-bound-genera2} and~\ref{cor:lower-bound-genera3} is strictly positive.  For instance,  the choice
$\upsilon=64(d+1)$  is sufficient and assumed throughout the proof.
\par
In what follows
$C_1,C_2,\ldots,$ and $c_1,c_2,\ldots,$ denote  constants that
may depend on $\vec{\beta}, \vec{r}, q, Q$ and $\Lambda$, but they are independent of $\vec{L}$ and $n$.
\subsection{Verification of Assumption \ref{ass1:mixed}}
Let $C:=\int_{-1}^1e^{-\frac{1}{1-z^2}}\rd z$, and
$$
U(x):=C^{-d}e^{-\sum_{j=1}^d\frac{1}{1-x_i^2}}\mathrm{1}_{[-1,1]^d}(x),\;\; x=(x_1,\ldots,x_d)\in\bR^d,
$$
For  $N>0$ and $a>0$ define
\begin{gather*}
\bar{f}_{0,N}(x):=(N)^{-d}\int_{\bR^d}U(y-x)\mathrm{1}_{[-N-1,-1]^d}(y)\rd y,\;\;\;
f_{0,N}(x):=a^d\bar{f}_{0,N}\big(xa\big).
\end{gather*}
\begin{lemma}
\label{lem:density-f_0,N}
The following statements hold.
\par
{\rm (a).} For any $N$ and $a$, $f_{0,N}$ is a probability density.
For any $\vec{\beta}, \vec{L}\in (0,\infty)^d$ and  $\vec{r}\in (0,\infty]^d$
there exists $a>0$ such that
$$
f_{0,N}\in\bN_{\vec{r},d}\big(\vec{\beta},\tfrac{1}{2}\vec{L}\big),\quad\forall N>0.
$$
\par
{\rm (b).} For any $Q>0$ and $q\in (r^*,\infty]$ there exists $N(q,Q)>0$ such that
$$
f_{0,N}\in\bN_{\vec{r},d}\big(\vec{\beta},\tfrac{1}{2}\vec{L}\big)\cap\bB_q(\tfrac{1}{2}Q),\quad\forall N\geq N(q,Q).
$$
\end{lemma}
\noindent The proof of the lemma is trivial; it is  omitted.
Let
\begin{equation}
\label{eq:def-f_0}
f_{0}:=f_{0,N},\;\quad  N\geq N(q, Q),
\end{equation}
so that statements (a) and (b) of Lemma~\ref{lem:density-f_0,N} hold.
\par
For  $\mu,\nu\in\mP[0,1]$ and  $z>0$ define
$e^*(z):=\max[e_\mu(z), e_\nu(z)]$.
\begin{lemma}
\label{lem:verification-ass-on-family}
Let $f_{0}$ be given in (\ref{eq:def-f_0}) and $f_w$ be defined in (\ref{eq:f-w})
 with
 $\Lambda\in\bC^{\infty}(\bR^d)$ satisfying~(\ref{eq:Lambda}),  and $\vec{\sigma}\in (0,1]^d$. Let  $\mu,\nu\in\mP[0,1]$ satisfy (\ref{eq:e>1/2}) and assume that
\begin{eqnarray}
\label{eq3*-lbf}
&&\mA\|\Lambda\|_q\big[\bma Me^*(q)\big]^{1/q}\leq  Q/4;
\\*[2mm]
\label{eq4-lbf}
&&\mA\sigma_l^{-\beta_l} \|\Lambda\|_{r_l}
\big[\bma Me^*(r_l)\big]^{1/r_l}\leq c_1  L_l,\quad\forall l=1,\ldots, d.
\end{eqnarray}
Then Assumption~\ref{ass1:mixed} is fulfilled for
$\bF=\bN_{\vec{r}, d}(\vec{\beta}, \vec{L})\cap \bB_q(Q)$ with  $\e=2^{-6}$.
\end{lemma}
The proof is given in Appendix.
\subsection{Key proposition}
The sought lower bound  depends
on parameters $A, \bma, M$ of family $\{f_w, w\in [0,1]^M\}$,
and parameters that characterize properties  of the prior measures
$\mu$ and $\nu$ in Definitions~\ref{def2:prior} and ~\ref{def3:prior}.
All these parameters should be specified to
satisfy conditions of Corollaries~\ref{cor:lower-bound-genera2} and~\ref{cor:lower-bound-genera3}.
It turns out that the choice of all parameters
can be made in a unified way so that the resulting lower bound in expressed only in
terms of the sample size $n$ and properties of measures $\mu$ and~$\nu$.
The corresponding statement is given in
Proposition~\ref{prop3} below; it is   of independent interest.
\par
For $r\in \bN^*$ we define
\[
 \bva:=\left\{\begin{array}{ll}
              1/n, & p\in \bN^*,\\
              r/n, & p\notin \bN^*,
             \end{array}
\right.
\;\;
\Bt:=\left\{\begin{array}{ll}
              p , & p\in \bN^*,\\
              \infty, & p\notin \bN^*,
             \end{array}
             \right.
             \]
and
\[
\md_{\mu, \nu}:=\frac{64\upsilon[e^*(p)]^2}{|e_\mu(p)-e_\nu(p)|^{2}},\;\;\;
\mn_{r,n}:=\left \{\begin{array}{ll}
n^{p/(p-1)}, & p\in \bN^*,
\\
e^{r}\wedge (n/r), & p\notin \bN^*.
\end{array}
\right.
\]
Define also
$$
\mJ_{r, n}(\mu, \nu):=\big[\md_{\mu, \nu}, \mn_{r,n}\big]\cap\Big\{x>0:\; x^{1-\frac{2}{\Bt}}\bva^2n\leq 1,\;\big([e^*(q)]^{1/q}\bva\big)^{\tau(q)}
x^{\frac{1}{q}-1+(1-\frac{1}{\Bt})\tau(q)}\leq 1\Big\},
$$
and
$$
\BM_{r,n}(\mu,\nu):=\sup_{M\in\mJ_{r,n}(\mu,\nu)}
M^{\big[\frac{1}{p}-\frac{1}{\Bt}+(1-\frac{1}{\Bt})(\frac{1}{\beta p}-\frac{1}{\omega})\big]\frac{1}{\tau(1)}}.
$$
\par
Several remarks on the choice of parameters $A, \bma$ and $M$
that clarify the above  definitions are in order.
In our parameter choice we treat condition (\ref{eq4-lbf}) of Lemma~\ref{lem:verification-ass-on-family},
the first condition in~(\ref{eq:condition-main2}) of Corollary~\ref{cor:lower-bound-genera2}
for  $p\notin\bN^*$, or condition (\ref{eq:condition-main})~of Corollary \ref{cor:lower-bound-genera3} for
$p\in\bN^*$ as {\em equalities}.
This allows us to  express  $A$ and $\sigma_l, l=1,\ldots d$ (and $\bma$)
as functions of $M, r$ and $n$. All other conditions
[Assumption~\ref{ass2:mixed}, condition (\ref{eq3*-lbf}) of
Lemma~\ref{lem:verification-ass-on-family} and the second bound in~(\ref{eq:condition-main2}) if $p\notin\bN^*$]
for given $r$ and $n$ determine the set to which parameter $M$ should belong.
In particular, the set
$$
\Big\{x>0:\; x^{1-2/\Bt}\bva^2n\leq 1,\;\big([e^*(q)]^{\frac{1}{q}}\bva\big)^{\frac{\tau(q)}{\tau(1)}}
x^{\frac{1/q-1+(1-1/\Bt)\tau(q)}{\tau(1)}}\leq 1\Big\}
$$
is determined by
conditions~(\ref{eq5-lbf}) of Assumption~\ref{ass2:mixed} and (\ref{eq3*-lbf})~of Lemma \ref{lem:verification-ass-on-family}. The quantity $\mn_{r, n}$
is related to condition (\ref{eq3-lbf}) of Assumption \ref{ass2:mixed} and to the second requirement
in  (\ref{eq:condition-main2}) of Corollary \ref{cor:lower-bound-genera2} if $p\notin\bN^*$.
It is readily seen that $\mJ_{r,n}(\mu,\nu)$ is the intersection of three intervals and, therefore, the value $\BM_{r,n}(\mu,\nu)$ is attained at
the one of their endpoints.  The quantity  $\md_{\mu,\nu}$ comes from the
condition $\Delta(\mu,\nu)>0$. We remark that  $\md_{\mu,\nu}>36\upsilon$ and, therefore, condition  (\ref{eq3**-lbf}) of Assumption \ref{ass2:mixed} is not active in the
$\bL_p$-norm estimation.

\begin{proposition}
\label{prop3}
Let  $r\in\bN^*$ be an integer number, possibly dependent on $n$,  and let $p>1$   be fixed.
Let  $\mu,\nu\in\mP[0,1]$ satisfy (\ref{eq:e>1/2}),
$\mJ_{r,p}(\mu,\nu)\neq\emptyset$, and  assume that
 $\ma\stackrel{p,p}{\backsim} \mb$ if $p\in\bN^*$, and
 $\ma\stackrel{2r}{\backsim} \mb$ if $p\notin\bN^*$. Then
\begin{align}
\cR_n\big[\bN_{\vec{r},d}\big(\vec{\beta},\vec{L}\big)&\cap\bB_q(Q)\big]
\nonumber
\\
&\geq C_1 \BL^{\frac{1-1/p}{\tau(1)}}\big(\bva\big[e^*(q)\big]^{\frac{1}{q}}\big)^{\frac{\tau(p)}{\tau(1)}}\BM_{r,n}(\mu,\nu)
\Bigg(\frac{|e_\mu(p)-e_\nu(p)|}{[e^*(q)]^{\frac{1}{q}}[e^*(p)]^{1-\frac{1}{p}}}\Bigg).
\label{eq:R>}
\end{align}
\end{proposition}
Proposition~\ref{prop3} relates properties of probability measures $\mu$ and $\nu$
as determined by
their moments and parameter $r$ with the sample sample $n$.
It provides a guideline for the choice of prior measures $\mu$ and $\nu$: they should be specified to
maximize the right hand side of (\ref{eq:R>}).

\subsection{Prior measure construction}\label{sec:prior-construction}
In view of Proposition \ref{prop3}, the proof of Theorems~\ref{th:lb-integer-r} and~\ref{th:lb-noninteger-r} is
reduced to the construction
of a pair of prior probability measures, $\mu, \nu\in \mP[0,1]$,
with required properties.
We present two statements, Propositions~\ref{prop1} and~\ref{prop2},
that state existence and provide explicit construction of such measures. Their proofs are given in Appendix.
\begin{proposition}
\label{prop1}
For any $t,s\in\bN^*$, $s\geq t>1$  one can construct  a pair of
probability measures $\mu, \nu\in \mP[0,1]$ such that  $e_\mu(z),e_\nu(z)\geq 1/2$ for any $z>0$,
$\mu\stackrel{s,t}{\backsim}\nu$ and
$$
e_{\mu}(t)-e_{\nu}(t)\geq  C_{s,t}:=\sqrt{2t-1}[(t-1)!]^2 (s-t)![(s+t-1)!]^{-1}.
$$
\end{proposition}
\par
For $S\in \bC(0,1)$ and $s\in \bN^*$  let
$\varpi_s(S)$ denote accuracy of the best approximation of $S$ on $[0,1]$ by algebraic
polynomials of degree $s$:
$$
\varpi_{s}(S):=\inf_{t\in\bR^{s+1}}\sup_{x\in[0,1]}\big|S(x)-P_{s,a}(x)\big|,
$$
where
$P_{s,a}(x):=\sum_{j=0}^{s}a_j x^j$, $a=(a_0, \ldots, a_s)\in \bR^{s+1}$.
\begin{proposition}
\label{prop2}
For any $S\in\bC(0,1)$ with $\varpi_s(S)>0$ there exist a pair of probability measures
$\mu, \nu\in \mP[0,1]$ such that  $e_\mu(z),e_\nu(z)\geq 1/2$ for any $z>0$,
$\mu\stackrel{s}{\backsim}\nu$ and
\[
\int_{0}^1 S(x) \mu(\rd x)-\int_{0}^1 S(x) \nu(\rd x)= \varpi_{s}(S).
\]
\end{proposition}
We remark that the measures $\mu,\nu$ constructed in Propositions \ref{prop1} and \ref{prop2} obviously satisfy the requirement (\ref{eq:e>1/2}).
\par
Lower and upper bounds for the accuracy of best approximation
$\varpi_{s}(S)$ are known for many continuous functions $S$.
In particular, the following results can be found in \cite{timan}, \S 7.1.41 and \S 7.5.4.
\begin{proposition}
\label{prop:timan}
For any  $p>0$,  $p\notin\bN^*$
there exists $C_p>0$ such that

\vskip0.15cm

\centerline{$
\varpi_{s}(x\mapsto x^p)\geq C_p s^{-2p},\quad\forall s\in\bN^*.
$}

\vskip0.15cm

 There exists $C>0$ such that

\vskip0.15cm

\centerline{
$
\varpi_{s}(x\mapsto x\ln x)\geq C s^{-2},\quad\forall s\in\bN^*.
$}
\end{proposition}
\subsection{Lower bounds corresponding to the regime $n^{-\frac{\tau(p)}{\tau(1)}}$}
In this case   lower bounds of Theorems~\ref{th:lb-integer-r} and~\ref{th:lb-noninteger-r} coincide. We derive them in a unified way below.

Put for brevity
$\me:=[e^*(q)]^{1/q}$ and assume that $\tau(q)\geq 0$.
\par
If $p\in \bN^*$  (under the premise of Theorem~\ref{th:lb-integer-r})
then we pick prior measures $\mu, \nu$
as in Proposition~\ref{prop1} with $s=t=p$.  If $p\notin \bN^*$
(under the premise of Theorem~\ref{th:lb-noninteger-r}), then we choose $\mu, \nu$
as in Proposition~\ref{prop2} with $s=2 s_0$, where $s_0$ is chosen from the relation
\begin{equation*}
s_0=\inf\big\{s\in\bN^*:\; \max\big[C^{-2}_{p,p},64\upsilon C_p^{-2}s^{4p}\big]\leq e^s\big\}.
\end{equation*}
Here $C_{p,p}$ and $C_p$ are the constants from Propositions \ref{prop1} and \ref{prop:timan} respectively.
\par
This choice guarantees that $0<\md_{\mu,\nu}< e^{s_0}$  and, therefore, choosing
$M=e^{s_0}$ we can assert that $M\in\mJ_{r,n}(\mu,\nu)$ for sufficiently large $n$. Indeed, if $\tau(q)>0$ we have
$$
\bva\to 0, \quad n\bva^2\to 0,\quad [\me\bva]^{\frac{\tau(q)}{\tau(1)}}\to 0,\;\; n\to\infty
$$
and, therefore $\mJ_{r,n}(\mu,\nu)\supseteq[\md_{\mu,\nu}, e^{s_0}]$ for $n$ large enough. By the same reason if $\tau(q)=0$
$$
\mJ_{r,n}(\mu,\nu)\supseteq\big[\md_{\mu,\nu}, e^{s_0}\big]\cap \big\{x\geq 1:\; x^{1-2/\Bt}\bva^2n\leq 1\big\}=\big[\md_{\mu,\nu}, e^{s_0}\big]
$$
for sufficiently large $n$.
It remains to note that  $\mu$ and $\nu$ constructed in Propositions \ref{prop1} and \ref{prop2} satisfy obviously (\ref{eq:e>1/2}) and $|e_\mu(p)-e_\nu(p)|\geq c_1$.
Thus,
applying  Proposition \ref{prop3} with $r=s_0$ we get
\begin{equation*}
\cR_n\big[\bN_{\vec{r},d}\big(\vec{\beta},\vec{L}\big)\cap\bB_q(Q)\big]\geq C_{1}\BL^{\frac{1-1/p}{\tau(1)}}
(1/n)^{\frac{\tau(p)}{\tau(1)}}.
\end{equation*}
This completes the proof of Theorems \ref{th:lb-integer-r} and \ref{th:lb-noninteger-r}
in the cases $\tau(q)\geq 0$, $\tau(p)\leq 1$ if $p\in\bN^*$ and $\tau(q)\geq 0$
$\tau(p)\leq 1-1/p$ if $p\notin\bN^*$.
\subsection{Proof of Theorem  \ref{th:lb-integer-r}}
Let   prior measures $\ma, \mb$ be   chosen according to
Proposition~\ref{prop1} with $s=t=p$. Hence $|e_\mu(p)-e_\nu(p)|\geq c_1$,  $\md_{\mu,\nu}=c_2$ and (\ref{eq:e>1/2}) is fulfilled. Remembering that
$\bva=n^{-1}$, $\Bt=p\geq 2$ we assert that
\begin{eqnarray*}
\mJ_{r,n}(\mu,\nu)&=&\Big[c_2, n^{\frac{p}{p-1}}\Big]\cap\Big\{x>0:\; [\me/n]^{\frac{\tau(q)}{\tau(1)}}
x^{\frac{1/q-1+(1-1/p)\tau(q)}{\tau(1)}}\leq 1\Big\}.
\end{eqnarray*}
In addition,  we deduce from Proposition \ref{prop3} that for any $M\in\mJ_{r,n}(\mu,\nu)$
\begin{equation}
\label{eq521-lbf}
\cR_n\big[\bN_{\vec{r},d}\big(\vec{\beta},\vec{L}\big)\cap\bB_q(Q)\big]\geq C_{2}\BL^{\frac{1-1/p}{\tau(1)}}
\bva^{\frac{\tau(p)}{\tau(1)}}M^{\frac{(1-1/p)(1/(\beta p)-1/\omega)}{\tau(1)}}.
\end{equation}
Now we consider two cases.
\par
(a).~
Assume first that $1/(\beta p)-1/\omega\geq 0$ that is equivalent to $\tau(p)\geq 1$.
Let us show that
$M:=n^{\frac{1}{1-1/p}}\in\mJ_{r,n}(\mu,\nu)$. Indeed,
$$
n^{-\frac{\tau(q)}{\tau(1)}}
M^{\frac{1/q-1+(1-1/p)\tau(q)}{\tau(1)}}=n^{\frac{1/q-1}{(1-1/p)\tau(1)}}\leq 1
$$
because $q\geq 1$. Thus,  we conclude from (\ref{eq521-lbf})
\begin{equation*}
\cR_n\big[\bN_{\vec{r},d}\big(\vec{\beta},\vec{L}\big)\cap\bB_q(Q)\big]\geq
C_{3}\BL^{\frac{1-1/p}{\tau(1)}}
(1/n)^{\frac{1}{\tau(1)}}.
\end{equation*}
\par
(b). Now let us assume that  $1/(\beta p)-1/\omega< 0$ and $\tau(q)< 0$. Let us show that
$$
M:=[\me n]^{\frac{\tau(q)}{(1-1/p)\tau(q)-(1-1/q)}}\in\mJ_{r,n}(\mu,\nu).
$$
Indeed, we have
$$
\frac{\tau(q)}{(1-1/p)\tau(q)-(1-1/q)}-\frac{1}{1-1/p}=\frac{(1-1/q)(1-1/p)^{-1}}{(1-1/p)\tau(q)-(1-1/q)}<0
$$
and, therefore, $M<n^{\frac{p}{p-1}}$. Moreover $M\to\infty, n\to\infty$.
Thus, $M\in\mJ_{r,n}(\mu,\nu)$  and we conclude  from (\ref{eq521-lbf})  that
\begin{equation*}
\cR_n\big[\bN_{\vec{r},d}\big(\vec{\beta},\vec{L}\big)\cap\bB_q(Q)\big]\geq C_{4}\BL^{\frac{1-1/p}{\tau(1)}}
(1/n)^{\frac{1/p-1/q}{1-1/q-(1-1/p)\tau(q)}}.
\end{equation*}
The theorem is proved.
\epr

\subsection{Proof of Theorem \ref{th:lb-noninteger-r}} Let  prior measures  $\ma, \mb$ be
chosen according to  Proposition~\ref{prop2} with $s=2\lfloor\ln n\rfloor+2$.  In Proposition \ref{prop3} choose $r=\lfloor\ln n\rfloor+1$;
then by Proposition~\ref{prop:timan}
\[
 |e_\mu(p)-e_\nu(p)|\geq C_1 r^{-2p}=c_3 (\ln n)^{-2p},\quad \md_{\mu,\nu}\leq c_4  (\ln n)^{4p}.
\]
Also, $\bva=\ln(n)/n$ and $\Bt=\infty$. Hence, we have for all $n$ large enough
\begin{equation*}
\mJ_{r,n}(\mu,\nu)=\Big[c_4(\ln n)^{4p}, n/\ln^{2}(n)\Big]\cap\Big\{x>0:\; [\me\ln(n)/n]^{\frac{\tau(q)}{\tau(1)}}
x^{\frac{1/q-1+\tau(q)}{\tau(1)}}\leq 1\Big\},
\end{equation*}
where we remind that $\me=[e^*(q)]^{1/q}$.
Since
$\mu$ and $\nu$ satisfy (\ref{eq:e>1/2})
 we derive from
Proposition \ref{prop3} that for any $M\in\mJ_{r,n}(\mu,\nu)$
\begin{equation}
\label{eq100:lbf-arbitrary-p}
\cR_n\big[\bN_{\vec{r},d}\big(\vec{\beta},\vec{L}\big)\cap\bB_q(Q)\big]\geq C_5\BL^{\frac{1-1/p}{\tau(1)}}(\ln(n)/n)^{\frac{\tau(p)}{\tau(1)}}M^{\frac{1/p+1/(\beta p)-1/\omega}{\tau(1)}}(\ln n)^{-2p}.
\end{equation}
Consider now two cases. First we assume first that $1/p+1/(\beta p)-1/\omega\geq 0$ (which
is equivalent to $\tau(p)\geq 1-1/p$) and show that
$M:=n/\ln^2(n)\in\mJ_{r,n}(\mu,\nu)$. Indeed,
$$
n^{-\frac{\tau(q)}{\tau(1)}}
M^{\frac{1/q-1+\tau(q)}{\tau(1)}}=\big[n/\ln^2(n)\big]^{\frac{1/q-1}{\tau(1)}} [\ln(n)]^{\frac{2-2/q}{\tau(1)}}\leq 1
$$
for all $n$ large enough. Hence, $M\in\mJ_{r,n}(\mu,\nu)$ and we conclude from
(\ref{eq100:lbf-arbitrary-p}) that
\begin{equation*}
\cR_n\big[\bN_{\vec{r},d}\big(\vec{\beta},\vec{L}\big)\cap\bB_q(Q)\big]\geq C_{6}\BL^{\frac{1-1/p}{\tau(1)}}
(1/n)^{\frac{1-1/p}{\tau(1)}}(\ln n)^{\frac{2(1-1/p)-\tau(p)}{\tau(1)}-2p}.
\end{equation*}
Second,
let $\tau(q)<0$ and show that
$
M:=\big(n/[\me\ln(n)]\big)^{-\frac{\tau(q)}{1-1/q-\tau(q)}}\in\mJ_{r,n}(\mu,\nu).
$
 Indeed,
$$
-\frac{\tau(q)}{1-1/q-\tau(q)}=1-\frac{1-1/q}{1-1/q-\tau(q)}\in (0,1),
$$
and, therefore, $M\in \big[c_4(\ln n)^{4p}, n/\ln^{2}(n)\big]$ for all $n$ large enough.
Hence $M\in\mJ_{r,n}(\mu,\nu)$
and we deduce from (\ref{eq100:lbf-arbitrary-p})
\begin{equation*}
\label{eq7:lbf-arbitrary-p}
\cR_n\big[\bN_{\vec{r},d}\big(\vec{\beta},\vec{L}\big)\cap\bB_q(Q)\big]\geq C_{7}\BL^{\frac{1-1/p}{\tau(1)}}
(1/n)^{\frac{1/p-1/q}{1-1/q-\tau(q)}}[\ln(n)]^{\frac{1/p-1/q}{1-1/q-\tau(q)}-2p}.
\end{equation*}
This completes the proof of the theorem.
\epr
\subsection{Proof of Proposition \ref{prop3}} The proof consists of  two steps.

\paragraph{Computation of  $\Delta(\ma,\mb)$}  In the specific
case of $\Psi(f)=\|f\|_p$ we have $G(x)=x^{1/p}$, $H(x)=x^p$ so that
functions $S_0$ and $S$ in (\ref{eq:S0-S}) and quantities
$E_\mc(A)$, $V_\pi(A)$ in (\ref{eq:E-pi-V-pi}) take the form
$S_0(z)=(1-z)^p \|f_0\|_p^p$, $S(z)= z^p \|\Lambda\|_p^p$, and, correspondingly,
\[
 E_\mc(A)= A^p \|\Lambda\|_p^p e_\mc(p),\;\;\; V_\mc (A)=A^p \|\Lambda\|_p^p [e_\mc(2p)]^{1/2}\leq 2A^p \|\Lambda\|_p^p e_\mc(p).
\]
The last inequality follows from (\ref{eq:e>1/2}). Therefore
\[
 H^*_{\mc}= [1-A\bma M e_\mc(1)]^p \|f_0\|_p^p + A^p\bma \|\Lambda\|_p^p M e_\mc(p).
\]
Since the choice of parameters $A$, $\bma$ and $M$   satisfies
(\ref{eq3-lbf}),
$1-A\bma M e_{\mc}(1)\geq 1/2$.
By definition of $f_0$, $\|f_0\|_p^p=c_2 N^{-d(p-1)}$, and $N$ can be chosen arbitrarily large.

In particular, denoting $e_{*}(p)=\min[e_\mu(p),e_\nu(p)]$ and choosing
 $N$ so that $N^{-d(p-1)}=c_3A^p\bma \|\Lambda\|_p^p e_{*}(p)\sqrt{\upsilon M}$
with sufficiently small $c_3>0$ we obtain
%
\begin{equation}\label{eq:H*}
A^p\bma \|\Lambda\|_p^p M e_\mc(p)\leq  H^*_{\mc} \leq A^p\bma \|\Lambda\|_p^p M e_\mc(p)+
A^p\bma \|\Lambda\|_p^p e_{*}(p)\sqrt{\upsilon M}.
\end{equation}
Furthermore,
\[
 |S_0(z)-S_0(z^\prime)|= \|f_0\|_p^p\, \big|(1-z)^p - (1-z^\prime)^p|\leq p \|f_0\|_p^p \;\big|z-z^\prime
 \big|,\;\;\;\forall z, z^\prime \in [0,1],
\]
so that
\[
 \eta_{S_0}\big(A\bma M\,e_\mc(1); A\bma \sqrt{\upsilon M}\,\big) \leq
 p \|f_0\|_p^p A\bma \sqrt{\upsilon M} \leq c_4N^{-d(p-1)}A\bma \sqrt{\upsilon M}.
\]
Taking into account that $N^{-d(p-1)}=c_3A^p\bma \|\Lambda\|_p^p e_{*}(p)\sqrt{\upsilon M}$ we obtain
\begin{eqnarray*}
 \alpha_\pi &\leq& c_4 A^{p+1}\bma^2 \|\Lambda\|_p^p e_{*}(p)\upsilon M + 2A^p \bma\|\Lambda\|_p^p e_\pi(p)\sqrt{\upsilon M}
\\*[2mm]
&\leq& 2A^p \bma\|\Lambda\|_p^p e_\pi(p)\sqrt{\upsilon M} \big(1+ (c_4/2)A\bma \sqrt{\upsilon M}\big)\leq
3 A^p\bma \|\Lambda\|_p^p e_\pi(p)\sqrt{\upsilon M},
\end{eqnarray*}
where in the last inequality we have used
(\ref{eq5-lbf}), and $\kappa_0$ is sufficiently small.
Therefore (\ref{eq:H*}) and (\ref{eq3**-lbf}) imply
\begin{align*}
H_{\mc}^* - \alpha_\pi&\geq A^p\bma \|\Lambda\|_p^p M e_\mc(p) \big(1-3\sqrt{\upsilon/M}\big),
\\
H_{\mc}^* + \alpha_\pi &\leq  A^p\bma \|\Lambda\|_p^p M e_\mc(p) \big(1+4\sqrt{\upsilon/M}\big).
\end{align*}
Since, $G(x)=x^{1/p}$ we can assert that
\[
 \cJ_{\mc}\subseteq \Big[ A (M\bma)^{1/p} \|\Lambda\|_p  [e_\mc(p)]^{1/p}\big(1-3\sqrt{\upsilon/M}\big)^{1/p},
 A (M\bma)^{1/p} \|\Lambda\|_p  [e_\mc(p)]^{1/p}\big(1+4\sqrt{\upsilon/M}\big)^{1/p}\Big]
\]
and, therefore, denoting $(y)_+=\max[0,y]$, we get
$$
\Delta(\mu,\nu)\geq A (M\bma)^{1/p} \|\Lambda\|_p\Big([e^*(p)]^{1/p}\big(1-3\sqrt{\upsilon/M}\big)^{1/p}-
[e_*(p)]^{1/p}\big(1+4\sqrt{\upsilon/M}\big)^{1/p}\Big)_{+}.
$$
Assuming   that
\begin{equation}
\label{eq:Delta>0}
 |e_\mu(p)-e_\nu(p)|>8 e^*(p) \sqrt{\upsilon/M}
\end{equation}
we can guarantee that $\Delta(\mu,\nu)>0$ because intervals
$[H^*_\mu-\alpha_\mu, H^*_\mu+\alpha_\mu]$ and $[H^*_\nu-\alpha_\nu, H_\nu+\alpha_\nu]$ are
disjoint. Additionally, by the elementary inequality
\[
 \tfrac{1}{p}(a\vee b)^{-1+1/p}  |a-b| \leq |a^{1/p}-b^{1/p}|
 \;\;\;\forall a, b>0,\;p\geq 1
\]
applied with $a= e^*(p)\big(1-3\sqrt{\upsilon/M}\big)$ and
$b=e_*(p)\big(1+4\sqrt{\upsilon/M}\big)$ we get
$$
\Delta(\mu,\nu)\geq c_3 A (M\bma)^{1/p}\big|e_\mu(p)-e_\nu(p)\big|[e^*(p)]^{1/p-1}.
$$
With this lower bound on $\Delta(\mu,\nu)$ applying either
Corollary~\ref{cor:lower-bound-genera2} or Corollary \ref{cor:lower-bound-genera3}
 we come to the following lower bound  on the minimax risk in terms of  parameters $A$, $\bma$ and $M$
 of the family $\{f_w\}$ and properties of the probability measures $\mu$ and $\nu$:
\begin{equation}
\label{eq52-lbf}
\cR_n\big[\bN_{\vec{r},d}\big(\vec{\beta},\vec{L}\big)\cap\bB_q(Q)\big]\geq c_4 \mA (\bma M)^{1/p}\big|e_\mu(p)-e_\nu(p)\big|[e^*(p)]^{1/p-1}.
\end{equation}

\paragraph{Generic choice of parameters}\label{sec:generic-choice}
In this part we present the choice of  parameters $A$ and $\vec{\sigma}$ as functions of $M, r$, $n$ and
satisfying  conditions (\ref{eq5-lbf})--(\ref{eq3-lbf}) and (\ref{eq3*-lbf})--(\ref{eq4-lbf}).
\par
 For any $\Bt\in (1,\infty]$ and $\bva<1$ we let
\begin{eqnarray*}
\sigma_l&=&c^{1/\beta_l}_3
L_l^{-1/\beta_l}\BL^{\frac{1/\beta_l-1/(\beta_l r_l)}{\tau(1)}}
[\me\bva]^{\frac{\tau(r_l)}{\beta_l\tau(1)}}M^{\frac{1/r_l-1+(1-1/\Bt)\tau(r_l)}{\beta_l\tau(1)}};
\\*[1mm]
A&=&c_4\me^{-1}\BL^{\frac{1}{\tau(1)}}[\me\bva]^{\frac{1-1/\omega}{\tau(1)}}M^{-\frac{1/\Bt+(1-1/\Bt)/\omega}{\tau(1)}}
\\
\bma&=&c_3^{1/\beta}\BL^{-\frac{1}{\tau(1)}}[\me\bva]^{\frac{1/\beta}{\tau(1)}}M^{\frac{1/\omega-1/(\beta \Bt)}{\tau(1)}}.
\end{eqnarray*}
Simple algebra shows that (\ref{eq4-lbf}) holds if $c_4c_3^{1/(\beta r^*)-1}\leq c_2$ because
$[e^*(r_l)]^{1/r_l}\leq \me$ for all $l=1,\ldots, d$, and $q\geq r^*$. Moreover, it yields
\begin{gather}
\label{eq2223::proof-t1}
\;\; A\bma =c_4c_3^{1/\beta}\bva M^{-1/\Bt};
\\*[2mm]
\label{eq2001::proof-t1}
A\bma M=c_4c_3^{1/\beta}\bva M^{1-1/\Bt};
\\*[2mm]
\label{eq2001-new::proof-t1}
nA^2\bma^2M=c_4^2c_3^{2/\beta}
M^{1-2/\Bt}\bva^2n;
\\*[2mm]
\label{eq2002::proof-t1}
\qquad\qquad\qquad\quad A\me[\bma M]^{\frac{1}{q}}=c_4c_3^{1/(\beta q)}\BL^{\frac{1-1/q}{\tau(1)}}[\me\bva]^{\frac{\tau(q)}{\tau(1)}}
M^{\frac{1/q-1+(1-1/\Bt)\tau(q)}{\tau(1)}}.
\end{gather}
In addition,
\begin{eqnarray}
\label{eq3:proof-t1}
&& A (\bma M)^{1/p}=c_5\me^{-1}\BL^{\frac{1-1/p}{\tau(1)}}[\me\bva]^{\frac{\tau(p)}{\tau(1)}}
M^{\frac{1}{\tau(1)}\big[\frac{1}{p}-\frac{1}{\Bt}+(1-\frac{1}{\Bt})(\frac{1}{\beta p}-\frac{1}{\omega})\big]}.
\end{eqnarray}
 Introduce
$$
\mX_{\bva,\Bt}=\Big[36\upsilon, \bva^{-\frac{t}{t-1}}\Big]\cap\Big\{x>0:\; x^{1-2/\Bt}\bva^2n\leq 1,\;[\me\bva]^{\frac{\tau(q)}{\tau(1)}}
x^{\frac{1/q-1+(1-1/\Bt)\tau(q)}{\tau(1)}}\leq 1\Big\}.
$$
First of all we assert that in view of (\ref{eq2001::proof-t1})--(\ref{eq2002::proof-t1})  $M\in\mX_{\bva,\Bt}$ implies the verification of
Assumption \ref{ass2:mixed} and (\ref{eq3*-lbf}) if one chooses $c_4$ sufficiently small.

Our goal now is to show that for any $l=1,\ldots,d$
\begin{gather}
\label{eq2222::proof-t1}
\sigma_l\leq 1,\quad \forall M\in\mX_{\bva,\Bt}.
\end{gather}
Put $T_l=c^{1/\beta_l}_3
L_l^{-1/\beta_l}\BL^{\frac{1/\beta_l-1/(\beta_l r_l)}{\tau(1)}}$ and consider separately two cases.

Let $\tau(q)\geq 0$. In this case $\tau(r_l)\geq 0$ for any $l=1,\ldots, d$ since $q\geq r^*$ and we have
\begin{eqnarray*}
&&\sigma_l=T_l
\Big[\me\bva M^{1-1/\Bt}\Big]^{\frac{\tau(r_l)}{\beta_l\tau(1)}}M^{\frac{1/r_l-1}{\beta_l\tau(1)}}\leq \Big[\bva M^{1-1/\Bt}\Big]^{\frac{\tau(r_l)}{\beta_l\tau(1)}}M^{\frac{1/r_l-1}{\beta_l\tau(1)}}\leq T_l\leq 1,
\end{eqnarray*}
if one chooses $c_3$ sufficiently small.
Here we have used also that $r_l\geq 1$, $\me\leq 1$, $M\geq 1$ and $M\in\mX_{\bva,\Bt}$.

Let now $\tau(q)<0$. First we note that in this case  necessarily
$$
\me\bva\leq M^{-\frac{1/q-1+(1-1/\Bt)\tau(q)}{\tau(q)}}\;\;\Leftrightarrow\;\;\me\bva M^{1-1/t}\leq M^{-\frac{1/q-1}{\tau(q)}}.
$$
On the other hand
\begin{eqnarray*}
\sigma_l&=&T_l[\me\bva]^{\frac{\tau(q)}{\beta_l\tau(1)}}[\me\bva]^{\frac{\tau(r_l)-\tau(q)}{\beta_l\tau(1)}}
M^{\frac{1/r_l-1/q+(1-1/t)(\tau(r_l)-\tau(q))}{\beta_l\tau(1)}}
M^{\frac{1/q-1+(1-1/\Bt)\tau(q)}{\beta_l\tau(1)}}
\\
&=&T_l[\me\bva]^{\frac{\tau(q)}{\beta_l\tau(1)}}M^{\frac{1/q-1+(1-1/\Bt)\tau(q)}{\beta_l\tau(1)}}\Big[\me\bva M^{1-1/t}\Big]^{\frac{\tau(r_l)-\tau(q)}{\beta_l\tau(1)}}M^{\frac{1/r_l-1/q}{\beta_l\tau(1)}}
\\
&\leq& T_l M^{\frac{1/r_l-1/q}{\beta_l\tau(1)}\big[1-\frac{1/q-1}{\beta\tau(q)}\big]}=T_l M^{\frac{1/r_l-1/q}{\beta_l\tau(q)}}\leq T_l\leq 1,
\end{eqnarray*}
if one chooses $c_3$ sufficiently small.  Here we have used once again that $q\geq r_l$ for any $l=1,\ldots,d$, $M>1$ and $\tau(q)<0$.
Thus, (\ref{eq2222::proof-t1}) is proved, and
all assumptions of Theorem~\ref{th:lower-bound-general} and
Lemma~\ref{lem:verification-ass-on-family} are verified if $M\in\mX_{\bva,\Bt}$.
\par
Note that if
$\Bt=p\mathrm{1}_{\bN^*}(p)+\infty\mathrm{1}_{\bar{\bN}^*}(p)$ and $\bva=n^{-1}\mathrm{1}_{\bN^*}(p)+ rn^{-1}\mathrm{1}_{\bar{\bN}^*}(p)$ then
$$
\mJ_{r,n}(\mu,\nu)\subset\mX_{\bva,\Bt}.
$$
To get the latter inclusion we have used also that (\ref{eq:Delta>0}) is equivalent to $M\geq \md_{\mu,\nu}$ and $\md_{\mu,\nu}\geq 36\upsilon$.
Moreover, we deduce from (\ref{eq2223::proof-t1}) that assumptions (\ref{eq:condition-main2}) of Corollary \ref{cor:lower-bound-genera2} and
(\ref{eq:condition-main}) of Corollary \ref{cor:lower-bound-genera3} with $t=p$ are verified for sufficiently small $c_4$. The assertion of the proposition follows now from (\ref{eq52-lbf})  and (\ref{eq3:proof-t1}).
\epr

\section{Proofs of Theorem \ref{th:lower-bound-general} and Corollaries \ref{cor:lower-bound-genera2}--\ref{cor:lower-bound-genera3}}
\label{sec:proofs_Tms-general}
\subsection{Proof of Theorem \ref{th:lower-bound-general}}
\label{sec:subsec-proofs_T1}
We break the proof into several steps.
\subsubsection*{$1^0$. Product form bounds for $p_{\bze}(x)$}
Our first  step is to develop tight bounds on $\mmp_{\bze}(x)$ for all
$x\in\bR^{dn}$  possessing a product form structure
with respect to the coordinates of $\bze$.
Recall that $p_{\bze}(x)=\prod_{i=1} f_{\bze}(x_i)$ where $f_{\bze}$ is
defined in (\ref{eq:f-w}).
\par
Let
$\Lambda_{0}(\cdot):=f_{0}(\cdot)\big[1- \mA\bma\varrho_{\bze}(1)\big]$.
As it was mentioned above,  (\ref{eq3-lbf}) implies
$1- \mA\bma\varrho_{\bze}(1)>0$.
Because $\Pi_0\cap \Pi_m=\emptyset$ for all $m\in \cM$, and $\Pi_j\cap \Pi_m= \emptyset$ for all $m, j\in \cM$, $m\ne j$
we have the following representation of function $f_\bze(\cdot)$: for any $y\in \bR^{d}$
\begin{eqnarray*}
f_{\bze}(y)&=&\Lambda_{0}(y)\mathrm{1}_{\Pi_{0}}(y)+
\mA\sum_{m\in\cM}\bze_{m}\Lambda_{m}(y)\mathrm{1}_{\Pi_{m}}(y)
\\
&=&
\Lambda_{0}(y)^{1_{\Pi_{0}}(y)}\prod_{m\in\cM}
\big[\mA \bze_{m}\Lambda_{m}(y)\big]^{\mathrm{1}_{\Pi_{m}}(y)}.
\end{eqnarray*}
Therefore
\begin{gather*}
 \mmp_{\bze}(x)=
 \prod_{i=1}^n f_{\bze}(x_i)=
 \Big[\prod_{i=1}^n\Lambda_{0}(x_i)^{\mathrm{1}_{\Pi_{0}}(x_i)}\Big]\Big[\prod_{m\in\cM}
 \prod_{i=1}^n\big[\mA \bze_{m}\Lambda_{m}(x_i)\big]^{1_{\Pi_{m}}(x_i)}
\Big].
\end{gather*}
If we put
\begin{equation}\label{eq:T}
T(x):=\prod_{i=1}^n\Big\{\big[f_{0}(x_i)\big]^{\mathrm{1}_{\Pi_{0}}(x_i)}\prod_{m\in\cM}
\big[\mA\Lambda_{m}(x_i)\big]^{\mathrm{1}_{\Pi_{m}}(x_i)}
\Big\}
\end{equation}
then
\begin{equation}\label{eq:pw}
\mmp_{\bze}(x)=T(x)
\big[1- \mA\bma\varrho_{\bze}(1)\big]^{n_0(x)}\prod_{m\in\cM}\big[\bze_{m}\big]^{n_{m}(x)}.
\end{equation}
Our current goal is to derive bounds on $1- \mA\bma\varrho_{\bze}(1)$.
Recall that $\bE_{\mc} \rho_{\bze}(1)=Me_\mc(1)$, and
denote for brevity
\[
b:=\mA\bma,  \;\;u:=\bE_{\mc} \varrho_{\bze}(1)=M e_\mc(1),\;\;
D: = b/(1-bu),
\]
where $D$ is defined in (\ref{eq:bayes-ration-second2}).
For $\pi \in \{\mu, \nu\}$  define the set
\begin{eqnarray}\label{eq:W-pi}
\cW_{\mc}:=\Big\{w\in[0,1]^{M}:\; \big|\varrho_{w}(1)- M
e_{\mc}(1)\big|\leq \sqrt{\upsilon M}\Big\},
\end{eqnarray}
and  suppose that $\bze \in \cW_{\mc}$.
Then (\ref{eq3-lbf}) implies
\begin{equation}\label{eq:rho-u}
|\varrho_\bze(1) - u| \leq \sqrt{\upsilon M} \leq M,\;\;\;\;b|\varrho_{\bze}(1) - u|\leq 1/2.
\end{equation}
Moreover, by (\ref{eq3-lbf}) and $e_{\mc}(1)\leq 1$,
$bu= A\bma Me_{\mc}(1)\leq 1/2$, and  $D\leq 2b$.
We have
\begin{align}\label{eq:rho-1}
 1- b\varrho_{\bze}(1)= 1- bu -b\big[\varrho_{\bze}(1)-u\big]=(1-bu)\big(1- D[\varrho_{\bze}(1)-u]\big).
\end{align}
Using elementary inequality $1-t\leq e^{-t}$ we obtain from (\ref{eq:rho-1})
\begin{equation}\label{eq:pw-u}
 1-b\varrho_{\bze}(1)\leq (1-bu)e^{Du} \exp\{- D\varrho_{\bze}(1)\}.
\end{equation}
On the other hand,
taking into account that $D|\varrho_{\bze}(1)-u|\leq 1$ on the event $\{\bze\in \cW_{\mc}\}$ and applying
inequality $1-t\geq e^{-t}-\frac{1}{2}t^2$, $\forall t\geq -1$ we get
\begin{align*}
 1- D[\varrho_{\bze}(1)&-u] \geq  e^{Du} \exp\{-D\varrho_{\bze}(1)\}\Big(1- \tfrac{1}{2}D^2[\varrho_{\bze}(1)-u]^2e^{D[\varrho_{\bze}(1)-u]}\Big)
\\
&\geq  e^{Du} \exp\{-D\varrho_{\bze}(1)\}\Big(1- 2b^2[\varrho_{\bze}(1)-u]^2e^{2b[\varrho_{\bze}(1)-u]}\Big)
\\
&\geq  e^{Du} \exp\{-D\varrho_{\bze}(1)\}\big(1- 2eb^2\upsilon M \big)
\geq e^{Du} \exp\{-D\varrho_{\bze}(1)\}\big(1- 1/n\big),
 \end{align*}
where the second inequality follows from $D\leq 2b$, the third inequality is a consequence of  (\ref{eq:rho-u}),
and the last inequality follows from condition
(\ref{eq5-lbf}) with small enough $\kappa_0$ satisfying $2e\kappa_0^2 \upsilon \leq 1$.
This together with (\ref{eq:rho-1}) yields
\begin{equation}\label{eq:pw-l}
 1-b\varrho_{\bze}(1) \geq  (1-bu)e^{Du} \exp\{-D\varrho_{\bze}(1)\}\big(1- 1/n\big).
\end{equation}
Combining (\ref{eq:pw-u}) and (\ref{eq:pw-l}) with
(\ref{eq:pw}) and  $n_0(x)\leq n$ we get
$e^{-1} p_{\bze}^*(x) \leq p_{\bze}(x) \leq p_{\bze}^*(x)$, $\forall x\in \bR^{dn}$,
where
\begin{align}\label{eq:666AUX}
 p_{\bze}^*(x)&:=T(x) (1-bu)^{n_0(x)} e^{Du n_0(x)}
 \prod_{m\in \cM} e^{-Dn_0(x) \bze_{m}} \big[\bze_{m}\big]^{n_m(x)}.
\end{align}
Thus we showed that under conditions (\ref{eq5-lbf}) and (\ref{eq3-lbf})
\begin{equation}
\label{eq1:bounds-for-bxi}
 \big\{\bze\in \cW_{\mc}\big\} \subseteq \big\{e^{-1} p_{\bze}^*(x) \leq p_{\bze}(x) \leq p_{\bze}^*(x)\big\},\;\;\;
 \forall x\in \bR^{dn}.
\end{equation}
Since $\bze_{m}, m\in\cM$ are independent random variables we get
$$
\bE_{\mc} \big\{\mmp^*_{\bze}(x)\big\}=T(x)(1-bu)^{n_0(x)} e^{Du n_0(x)}
\prod_{m\in\cM}\bE_{\mc} \Big\{e^{-D n_{0}(x) \bze_{m}}\big[ \bze_{m}\big]^{n_{m}(x)}\Big\},
$$
It remains to note that since $\ma\backsim \mb$,  values of $D$ and $u$
do not depend on $\mc\in\{\ma, \mb\}$; therefore
\begin{eqnarray*}
\Upsilon(x):=
\frac{\bE_{\ma}\big\{ \mmp^*_{\bze}(x)\big\}}{\bE_{\mb}\big\{ \mmp^*_{\bze}(x)\big\}}=
\prod_{m\in\cM}\frac{\gamma_{m,\ma}(x)}{\gamma_{m,\mb}(x)},\quad\forall x\in\bR^{dn},
\end{eqnarray*}
where $\gamma_{m, \mc}$ is given in (\ref{eq:bayes-ration-second2}).
\subsubsection*{$2^0$. Derivation of lower bound (\ref{eq:lower-bound})}
(a).~According to (\ref{eq:int-H}),
\begin{eqnarray*}
\int_{\bR^d}H\big(f_{\bze}(x)\big)\rd x= S_0\big(A\bma\varrho_{\bze}(1)\big) +
\bma\sum_{m\in\cM}S\big(A\bze_{m}\big).
\end{eqnarray*}
For $\pi \in \{\mu, \nu\}$ define
\begin{eqnarray*}
\cW_{\mc, S}:= \Big\{w\in [0,1]^{M}: \Big|\sum_{m\in \cM} S(Aw_m) - ME_\mc (A)\Big|
\leq \sqrt{\upsilon M}\,  V_{\mc}(A)
\Big\}.
\end{eqnarray*}
It is worth noting that   events $\{\bze\in \cW_{\mc}\}$  [see (\ref{eq:W-pi})] and $\{\bze\in \cW_{\mc, S}\}$
control deviations of sums of independent random variables
from their expectations, where the thresholds on the right  hand side in the definitions of
$\cW_{\mc}$ and $\cW_{\mc, S}$ are the upper bounds on the
standard deviations of the sum  inflated by a factor $\sqrt{\upsilon}$.
This fact allows to assert that by Chebyshev's inequality
\[
\bP_{\mc}\{\bze\not\in \cW_{\mc}\}\leq 1/\upsilon, \;\;
 \bP_{\mc}\{\bze\not\in \cW_{\mc, S}\}\leq 1/\upsilon.
 \]
\par
Assume that $\bze\in \cW_{\mc}\cap \cW_{\mc, S}$; then
\begin{align*}
 \Big| \int_{\bR^d} &H\big(f_{\bze}(x)\big)\rd x  - H^*_{\mc}\Big|
 \\
 &\;\;\leq  \big|S_0\big(A\bma\varrho_{\bze}(1)\big)
 - S_0(A\bma M e_\mc (1))\big|
 +  \bma \Big|\sum_{m\in \cM} S(A\zeta_m) - M E_\mc (A)\Big|,
 \end{align*}
 where $H_\mc^*$ is defined in (\ref{eq:H-pi}).
 We have $A\bma |\varrho_\zeta(1)- M e_\mc(1)|\leq A\bma \sqrt{\upsilon M}$  in view of
 $\bze\in \cW_{\pi}$ and (\ref{eq5-lbf});
therefore
\[
 \big|S_0\big(A\bma\varrho_{\bze}(1)\big)
 - S_0(A\bma M e_\mc (1))\big| \leq \eta_{S_0}\big(A\bma Me_\mc(1);\,A\bma \sqrt{\upsilon M}\big).
\]
If  $\bze\in \cW_{\mc}\cap \cW_{\mc, S}$ is realized then
\begin{equation*}
 \Big| \int_{\bR^d} H\big(f_{\bze}(x)\big)\rd x  - H^*_{\mc}\Big|  \leq
 \eta_{S_0}\big(A\bma Me_\mc(1); A\bma \sqrt{\upsilon M}\big)  +
 \bma \sqrt{\upsilon M}\,V_\mc(A)= \alpha_{\mc},
\end{equation*}
where $\alpha_{\mc}$ is defined in (\ref{eq:alpha-pi}).
Therefore we have shown that
\begin{equation}\label{eq:zeta-Psi}
\{\bze\in \cW_{\mc}\} \cap \{\bze \in \cW_{\mc, S}\} \subseteq
\{\Psi(f_\bze) \in \cJ_{\mc}\},
\end{equation}
where $\cJ_{\mc}$ is defined in (\ref{eq:cJ-pi}).
\par\medskip
(b).
For $\pi\in \{\mu, \nu\}$
define  $\cC_\mc:=\{\bze\in \cW_{\mc}\}\cap \{\bze \in \cW_{\mc, S} \}\cap \{f_\bze \in \bF\}$.
For the sake of brevity
in the subsequent proof we write $\Delta:=\Delta(\ma, \mb)$.
For arbitrary estimator $\tilde{\Psi}$ of $\Psi(f)$
we have
\begin{align}
2\sup_{f\in \bF\cap \mF} \bP_f \big\{ |\tilde{\Psi}- \Psi(f)| \geq
\tfrac{1}{3}\Delta\big\}
 \;\geq
  \bE_{\ma}\Big[ {\rm 1}\{\cC_\ma\} \bP_{f_\bze} \big\{|\tilde{\Psi}- \Psi(f_\bze)|\geq
 \tfrac{1}{3}\Delta \big\}\Big]
 \nonumber
 \\ +
   \bE_{\mb}\Big[ {\rm 1}\{\cC_\mb\} \bP_{f_\bze} \big\{|\tilde{\Psi}- \Psi(f_\bze)|\geq
  \tfrac{1}{3}\Delta\big\}\Big].
\label{eq:000}
  \end{align}
 Let  $a_{\mc}:= \inf_{|\alpha|\leq \alpha_{\mc}} G(\bar{H}_{\mc}+\alpha)$,
  $b_{\mc}:= \sup_{|\alpha|\leq \alpha_{\mc}} G(\bar{H}_{\mc}+\alpha)$
  so that $\cJ_{\mc}=[a_{\mc}, b_{\mc}]$.
  Therefore letting
  $\cI_{\mc}(\Delta):=
 [ a_{\mc}-\tfrac{1}{3}\Delta,
 b_{\mc} + \tfrac{1}{3}\Delta]$ we have from (\ref{eq:zeta-Psi})
\[
 \cC_\mc \cap \Big\{\tilde{\Psi}-\Psi(f_{\bze})|\geq \tfrac{1}{3}\Delta\Big\}
 \supseteq \cC_\mc \cap \big\{\tilde{\Psi}\notin \cI_{\mc}(\Delta)\big\}.
\]
This implies
\begin{align}
J_{\mc}&:= \bE_{\mc} \Big[ {\rm 1}\{\cC_\mc\} \bP_{f_\bze} \big\{|\tilde{\Psi}
-  \Psi(f_\bze)|
\geq \tfrac{1}{3}
 \Delta \big\}\Big]
\geq
 \bE_{\mc}\Big[ {\rm 1}\{\cC_\mc\} \bP_{f_\bze} \big\{
 \tilde{\Psi}\notin \cI_{\mc}(\Delta)\big\}
 \Big]
\nonumber
 \\
 &=
 \bE_{\mc}\Big[ {\rm 1}\{\cC_\mc\} \int_{\bR^{dn}}
 {\rm 1}\big\{\tilde{\Psi}(x)
 \notin \cI_{\mc}(\Delta) \big\} \mmp_\bze (x) \rd x  \Big]
 \nonumber
 \\
 &\geq e^{-1} \bE_{\mc}\Big[ {\rm 1}\{\cC_\mc\} \int_{\bR^{dn}}
 {\rm 1}\big\{\tilde{\Psi}(x)
 \notin \cI_{\mc}(\Delta)\big\}
 \mmp^*_\bze (x) \rd x  \Big]
\nonumber
 \\
 &\geq   e^{-1}\int_{\bR^{dn}}
 {\rm 1}\big\{\tilde{\Psi}(x)
 \notin \cI_{\mc}(\Delta) \big\}
 \bE_{\mc} [ \mmp^*_\bze (x)] \rd x  -
 e^{-1}\bE_{\mc} \big[ {\rm 1}\{\bar{\cC}_\mc\}\, \mmp^*_\bze\big],
 \nonumber
 \end{align}
 where $\bar{\cC}_\mc$ is the event complementary to $\cC_\mc$, and $\mmp^*_\bze:=\int_{\bR^{dn}}
 \mmp_\bze^*(x) \rd x$. In the third line we have used
 that $\mmp_{\bze} (x) \geq e^{-1} \mmp^*_{\bze}(x)$ for all $x\in \bR^{dn}$
 on the event $\{\bze\in \cW_{\mc}\}$.
Note that  by Chebyshev's inequality and in view of Assumption~\ref{ass1:mixed}
 \[
  \bP_{\mc}\{ \bar{\cC}_\mc\} \leq \bP_{\mc}\{ \bze\not\in \cW_{\mc}\} +
  \bP_{\mc}\{\bze \not\in\cW_{\mc, S}\} +
  \bP_{\mc}\{\bze \not \in \bF\}\leq 2\upsilon^{-1} +\varepsilon.
 \]
 Then
by the Cauchy--Schwarz
inequality
\begin{equation}\label{eq:Remainder}
\bE_{\mc} \big[ {\rm 1}\{\bar{\cC}_\mc\}\, p^*_\bze\big]\leq
\sqrt{(2\upsilon^{-1}+\varepsilon})
\max_{\mc\in \{\ma, \mb\}}\big\{\bE_{\mc}(\mmp_{\bze}^*)^2\big\}^{1/2}=: R
\end{equation}
which leads~to
 \begin{equation}
 J_{\mc} \geq   e^{-1}
 \int_{\bR^{dn}}
 {\rm 1}\big\{\tilde{\Psi}(x)
 \notin \cI_{\mc}(\Delta) \big\}
 \bE_{\mc} [ \mmp^*_\bze (x)] \rd x - e^{-1} R.
 \label{eq:J-pi}
\end{equation}
Furthermore,  we note that
\begin{align}\label{eq:J-mu}
  J_{\ma} + e^{-1}R
  &\;\geq  e^{-1}\int_{\bR^{dn}}
 {\rm 1}\big\{\tilde{\Psi}(x)
 \notin \cI_{\ma}(\Delta)\big\}
 \Upsilon(x) \bE_{\mb} [ p^*_\bze (x)] \rd x
 \nonumber
\\
&\;\geq (2e)^{-1} \int_{\bR^{dn}}
 {\rm 1}\big\{
 \tilde{\Psi}(x)\notin \cI_{\ma}(\Delta)\big\}
 {\rm 1}\{\Upsilon(x)\geq \tfrac{1}{2}\}
 \bE_{\mb} [ p^*_\bze (x)] \rd x,
 \end{align}
and that
for all $x\in \bR^{dn}$
\begin{equation}
 {\rm 1}\big\{\tilde{\Psi}(x)\notin \cI_{\ma}(\Delta)\big\}
+{\rm 1}\big\{\tilde{\Psi}(x)\notin \cI_{\mb}(\Delta)\big\} \geq 1.
 \label{eq:sum-indicators}
\end{equation}
The last inequality is an immediate consequence of definition of $\cI_{\mc}(\Delta)$ and
the fact that $\Delta=\Delta(\ma, \mb)>0$.
Therefore combining (\ref{eq:sum-indicators}),
(\ref{eq:J-mu}), (\ref{eq:J-pi}) and (\ref{eq:000}) we obtain
\begin{align*}
 \tfrac{1}{2}(J_{\ma} + J_{\mb})
 &\geq
 (4e)^{-1} \int _{\bR^{dn}} {\rm 1}\{\Upsilon(x)\geq \tfrac{1}{2}\}\bE_{\mb} [ \mmp^*_\bze (x)] \rd x
 -e^{-1}R
 \\
 &\geq (4e)^{-1} \int _{\bR^{dn}} {\rm 1}\{\Upsilon(x)\geq \tfrac{1}{2}\}
 \bE_{\mb} \Big[{\rm 1}\{\bze\in \cW_{\mb}\}\mmp_\bze (x)\Big] \rd x -
 e^{-1}R
\\
 & \geq (4e)^{-1} \bE_{\mb} \Big[ {\rm 1}(\bze\in \cW_{\mb}) \bP_{f_\bze}\big\{\Upsilon(X^{(n)})
 \geq \tfrac{1}{2}\big\}
\Big] - e^{-1}R
\\
& \geq (4e)^{-1} \bE_{\mb} \Big[ \bP_{f_\bze}\big\{\Upsilon(X^{(n)})
 \geq \tfrac{1}{2}\big\}
\Big] - (4e\upsilon)^{-1}  - e^{-1}R,
 \end{align*}
 where in the second line we have used  that $\mmp_{\bze}^*(x) \geq \mmp_{\bze}(x)$ for all
 $x\in \bR^{dn}$ on the event $\{\bze\in \cW_{\mb}\}$, see (\ref{eq1:bounds-for-bxi}).
This together with (\ref{eq:000}) and Chebyshev's inequality implies that
\begin{align}\label{eq:LB-0}
 [\Delta (\ma, \mb)]^{-2}\,& \cR^2_n[\bF]
 \nonumber
 \\
 & \geq (36e)^{-1} \bE_{\mb} \Big[\bP_{f_\bze}\big\{\Upsilon(X^{(n)})\geq \tfrac{1}{2}\big\}
\Big] - (36 e \upsilon)^{-1} - (9e)^{-1}R.
\end{align}
\subsubsection*{$3^0$. Bounding the remainder  in (\ref{eq:LB-0})}
In order to complete the proof of the theorem, in view of (\ref{eq:Remainder}), it remains to show that
\begin{equation}
\label{eq:new6-lb}
\bE_{\mc}\big(\mmp^*_{\bze}\big)^2\leq 2.
\end{equation}
Indeed, if  (\ref{eq:new6-lb}) is established then the theorem statement follows from (\ref{eq:LB-0})
and (\ref{eq:Remainder}).
\par\medskip
(a).~First, we note that in view of (\ref{eq:T}) and by definition of $n_0(x)$ and $n_m(x)$
\begin{equation*}
T(x) (1-bu)^{n_0(x)} e^{Du n_0(x)}
=\prod_{i=1}^n \bigg\{\Big[e^{Du}(1-bu)f_0(x_i)\Big]^{\mathrm{1}_{\Pi_{0}}(x_i)}
\prod_{m\in\cM}
\big[\mA\Lambda_{m}(x_i)\big]^{\mathrm{1}_{\Pi_{m}}(x_i)}\bigg\},
\end{equation*}
and
\begin{align*}
 \prod_{m\in \cM}
 e^{-Dn_0(x) \bze_{m}}
 \big[\bze_{m}\big]^{n_m(x)}=
 \prod_{i=1}^d
 \big[e^{-D\varrho_{\bze}(1)}\big]^{{\rm 1}_{\Pi_0(x_i)}}\prod_{m\in \cM} [\bze_{m}]^{{\rm 1}_{\Pi_m(x_i)}}.
\end{align*}
Therefore
\begin{multline*}
 \mmp^*_{\bze}(x) = \prod_{i=1}^n \Big\{ \Big[f_0(x_i) (1-bu) e^{-D[\rho_{\bze}(1)-u]}\Big]^{{\rm 1}_{\Pi_0(x_i)}}
 \prod_{m\in \cM} \big(A\Lambda_m(x_i) \bze_{m}\big)^{{\rm 1}_{\Pi_m}(x_i)}\Big\}
 \\
 =\prod_{i=1}^n \Big\{ \Big[ {\rm 1}_{\bar{\Pi}_0}(x_i)+ f_0(x_i) (1-bu) e^{-D[\rho_{\bze}(1)-u]}
 \Big] \Big[{\rm 1}_{\Pi_0}(x_i)+ \sum_{m\in \cM} A \bze_{m} \Lambda_m(x_i)\Big]\Big\}
 \\
 =\prod_{i=1}^n \Big[ f_0(x_i) (1-bu) e^{-D[\rho_{\bze}(1)-u]} + \sum_{m\in \cM} A \bze_{m} \Lambda_m(x_i)\Big],
\end{multline*}
and taking into account that $\int \Lambda_m(x)\rd x=\bma$ and $b= A\bma$
we obtain
\[
 \mmp_{\bze}^*= \Big[(1-bu) e^{-D[\rho_{\bze}(1)-u]} + b\rho_{\bze}(1)\Big]^n= \Big[
e^{-\frac{b[\varrho_{\bze}(1)-u]}{1-bu}}(1-bu)+b\varrho_{\bze}(1)\Big]^n.
\]
Denote
$\chi:=\varrho_{\bze}(1)-\bE_{\mc}\big\{\varrho_{\bze}(1)\big\}$. Since
$\bE_{\mc}\big\{\varrho_{\bze}(1)\big\}=u$ we have
\begin{equation}
\label{eq:new2-lb}
\mmp^*_{\bze}=\big[
e^{-\frac{b\chi}{1-bu}}(1-bu)+bu +b\chi\big]^n.
\end{equation}
Note that $e_{\mc}(1)\leq 1$; hence
 $0<bu\leq 1/2$ in view of (\ref{eq3-lbf}). Also since $\varrho_{\bze}(1)$ is a positive random variable
 $$
 \frac{b(\varrho_{\bze}(1)-u)}{1-bu}\geq -\frac{bu}{1-bu}\geq -1.
 $$
Since $e^{-t}\leq 1-t+\frac{1}{2}t^2$ for all $t\geq -1$, 
 we have
$e^{-\frac{b\chi}{1-bu}}(1-bu)\leq 1-bu-b\chi+ 2b^2\chi^2$ which
together with (\ref{eq:new2-lb}) leads to
\begin{eqnarray}
\label{eq:new3-lb}
\mmp^*_{\bze}
\leq \big[1+2b^2\chi^2 \big]^n.
\end{eqnarray}
Now we bound second moment of the random variable on the right hand side of (\ref{eq:new3-lb}).
\par\medskip
(b). We have
\begin{eqnarray*}
\label{eq:new4-lb}
&&\bE_\mc  \big(\mmp^*_{\bze}\big)^2 \leq \bE_\mc\big\{e^{4nb^2\chi^2}\big\}\leq 1+
\int_1^\infty\bP_\mc\Big(|\chi|\geq (4nb^2)^{-1/2}\sqrt{\ln y}\Big)\rd y.
\end{eqnarray*}
Since $\chi$ is the sum of  $M$ centered i.i.d. random variables taking values in $[0,1]$, applying
the Hoeffding inequality we  obtain
\begin{eqnarray*}
\label{eq:new5-lb}
&&\bE_\mc\big(\mmp^*_{\bze}\big)^2 \leq 1+2\int_{1}^\infty e^{-(2Mnb^2)^{-1}\ln y}\rd y.
\end{eqnarray*}
Then
(\ref{eq5-lbf}) implies
$2Mnb^2=2n\mA^2\bma^2 M \leq 2\kappa_0^2$,
and
\begin{eqnarray*}
\bE_\mc\big(\mmp^*_{\bze}\big)^2\leq 1+2\int_{1}^\infty y^{-1/(2\kappa_0^{2})}\rd y=1+
4\kappa^2_0(1-2\kappa_0^2)^{-1}\leq 2
\end{eqnarray*}
for sufficiently small $\kappa_0$.
Thus, (\ref{eq:new6-lb})  is proved.
\epr

\subsection{Proof of Corollary \ref{cor:lower-bound-genera2}}
\label{sec:subsec-proofs_T2}
The following well known inequality on the tail of binomial random variable
[see, e.g., in \cite{Bouch-Lugosi-Massart}]
will be exploited in the proof of the theorem.
\begin{lemma}
\label{lem:lb-deviation}
Let $\xi$ be a binomial random variable with parameters $n$ and $p$. Then for $pn\leq z\leq n$ one has
\[
 \bP (\xi \geq z) \leq \Big(\frac{pn}{z}\Big)^z e^{z-pn}.
\]
\end{lemma}
%
%
\noindent First of all let us remark that the first condition in (\ref{eq:condition-main2}) implies
\begin{equation}
\label{eq14:AUX}
Dn\leq 2\kappa_1 r.
\end{equation}
Also without loss of generality we will assume that $e_\mu(r)\leq e_\nu(r)$, since  all our constructions and definitions  are "symmetric" w.r.t $\mu$ and $\nu$.
We have
for any $m\in\cM$
\begin{equation*}
\label{eq44-lbf}
\gamma_{m,\mc}(x)=\int_{0}^1 y^{n_{m}(x)}e^{-Dn_{0}(x)y}\mc(\rd y)
=\sum_{k=0}^\infty\frac{(-1)^k}{k!} \big[Dn_{0}(x)\big]^k  e_{\mc}\big(n_{m}(x)+k\big).
\nonumber
\end{equation*}
Let
$$
Y_{m}(x):=\sum_{k=0}^{r}\frac{(-1)^k}{k!}
\big[Dn_{0}(x)\big]^k \big[e_\ma(n_{m}(x)+k)-e_\mb(n_{m}(x)+k)\big].
$$
Taking into account that
$\mu\stackrel{2r}{\backsim}\nu$,
$e_\pi(k_1)\geq e_\pi(k_2)$ for  $k_1\leq k_2$,  and (\ref{eq14:AUX})
we obtain
\begin{align}
\big|\gamma_{m,\ma}(x)&-\gamma_{m,\mb}(x)\big| \leq
\big|Y_{m}(x)\big|+\big[e_\mu(r)+
e_\nu(r)\big]\sum_{k=r+1}^\infty\frac{(Dn)^k}{k!}
\nonumber
\\
&\leq \big|Y_{m}(x)\big|+
2e_\nu(r)
\sum_{k=r+1}^\infty\bigg[\frac{2e\kappa_1 r}{k}\bigg]^k
\leq \big|Y_{m}(x)\big|+
4(2e\kappa_1)^{r+1}e_\nu(r)
\label{eq45-lbf}
\end{align}
for small enough
$\kappa_1$.
On the other hand, by definition
%
\begin{equation}
\label{eq46-lbf}
\gamma_{m,\mc}(x)\geq e^{-Dn}e_{\mc}\big(n_{m}(x)\big).
\end{equation}
Now  define the random event
\[
\cA:=\cup_{m\in\cM}\cA_{m},\;\;\;
\cA_{m}:=\big\{n_{m}\big(X^{(n)}\big)\geq r\big\},
\]
and let $\bar{\cA}$ be the event complimentary to $\cA$.

Since $\ma\stackrel{2r}{\backsim}\mb$, on the  event $\bar{\cA}$ we have
$Y_{m}\big(X^{(n)}\big)=0$, $\forall m\in\cM$, and
$e_{\mb}\big(n_{m}(x)\big)\geq e_{\mb}(r)$.
Therefore
it follows from
from (\ref{eq45-lbf}) and (\ref{eq46-lbf}) that
\begin{eqnarray}
\label{eq47-lbf}
\Upsilon\big(X^{(n)}\big)&\geq& \prod_{m\in\cM}
\big[1-4e^{2\kappa_1 r}(2e\kappa_1)^{r}\big]= \prod_{m\in \cM} \big[1-4e^{c_0r}\big]
\end{eqnarray}
where $c_0=2\kappa_1+1+\ln(2\kappa_1)$.
In view of the second condition in (\ref{eq:condition-main2}),
$M\leq e^{r}$, and  (\ref{eq47-lbf}) implies that  if  $\bar{\cA}$ is realized then
\begin{equation*}
\label{eq48-lbf}
\Upsilon\big(X^{(n)}\big)\geq (1-4e^{c_0 r})^{M}\geq \inf_{z\geq 2}\big(1-z^{c_0}\big)^{z} \geq \tfrac{1}{2},
\end{equation*}
where the last inequality holds because $c_0=2\kappa_1+1+\ln (2\kappa_1)$
can be made negative and arbitrary small  by choice of  sufficiently small $\kappa_1$.
Thus
\begin{equation}
\label{eq49-lbf}
\bar{\cA}\subseteq \big\{\Upsilon\big(X^{(n)}\big)> \tfrac{1}{2}\big\},
\end{equation}
and we derive from (\ref{eq49-lbf})
\begin{equation}
\label{eq490000-lbf}
\bP_{f_{\bze}} \big\{\Upsilon\big(X^{(n)}\big)> 2^{-1}\big\}\geq \bP_{f_{\bze}}\big(\bar{\cA}\big)\geq 1-\sum_{m\in\cM}\bP_{f_{\bze}}\big(\cA_{m}\big).
\end{equation}
Our current goal is to bound from below the expression on the right hand side of the last formula.
\par
Note that for any $m\in \cM$ random variable
$\widehat{n}_m=n_m(X^{(n)})=\sum_{i=1}^n \mathrm{1}_{\Pi_m}(X_i)$ has binomial distribution with parameters
$n$ and
$$
p_m:=\bP_{f_{\bze}}\big(X_i\in\Pi_{m}\big)=\int_{\Pi_m}f_{\bze}(y)\rd y=\mA\bma\bze_{m}\leq \mA\bma.
$$
By the first condition in (\ref{eq:condition-main2}), $p_mn \leq A\bma n\leq \kappa_1r\leq r$ so that we can
apply Lemma~\ref{lem:lb-deviation}:
\[
 \bP_{f_\bze}(\cA_m)=\bP_{f_\bze}\{\widehat{n}_m\geq r\}\leq (e A\bma n r^{-1})^r \leq (e\kappa_1)^{r}.
\]
Therefore taking into account that $M\leq e^r$ we obtain
\begin{equation}
\label{eq51-lbf}
\bP_{f_{\bze}}\big(\bar{\cA}\big)\geq 1- M(e\kappa_1)^{r}\geq 1-(e^2\kappa_1)^r \geq \tfrac{1}{3},
\end{equation}
provided that $\kappa_1$ is small enough.
Then it follows  from (\ref{eq490000-lbf}) and (\ref{eq51-lbf})
that
$$
\bE_{\mb}\big[
\bP_{f_{\bze}}\big\{\Upsilon\big(X^{(n)}\big)\geq \tfrac{1}{2}\big\}\big]\geq
\tfrac{1}{3}.
$$
The statement of the corollary follows now from Theorem \ref{th:lower-bound-general}.
\epr

\subsection{Proof of Corollary \ref{cor:lower-bound-genera3}}
\label{sec:subsec-proofs_T3}
Let $r=t$, and $\mu, \nu\in \mP[0,1]$ satisfy $\mu\stackrel{t,t}{\backsim}\nu$.
Let us remark that   (\ref{eq:condition-main}) implies
\begin{equation}
\label{eq9:AUX}
Dn\leq 2A\bma n\leq  2\kappa_1M^{-1/t}.
\end{equation}
Also without loss of generality we will assume that $e_\mu(t)\leq e_\nu(t)$.
\par
1$^0$.  Define the random events
$$
\cD:=\bigcap_{m\in\cM}\big\{\widehat{n}_{m}\leq t-1\big\},\quad \cE:=
\bigcap_{k=0}^{t-1}\Big\{\eta_k\leq M^{1-\frac{k}{t}}\Big\},
$$
where
$\widehat{n}_{m}=n_{m}(X^{(n)})$, and we have put $\eta_k:=\sum_{m\in\cM}\mathrm{1}\{\widehat{n}_{m}=k\}$.
If  $\cD$ is realized then for any $m\in\cM$  we have
\begin{equation*}
\gamma_{m,\pi}\big(X^{(n)}\big)=\sum_{k=0}^{t-1}\mathrm{1}\big(\widehat{n}_{m}=k\big)
\bigg[\int_{0}^1 y^ke^{-\widehat{n}_{0} D y}\pi(\rd y)\bigg]=\prod_{k=0}^{t-1}\big[\mT_k(\pi)\big]^{\mathrm{1}(\widehat{n}_{m}=k)},
\end{equation*}
where, recall, $\gamma_{m,\pi}(\cdot)$ is defined in (\ref{eq:bayes-ration-second2}), and we have denoted
$$
\mT_k(\pi):=\int_{0}^1 y^ke^{-\widehat{n}_{0} D y}\pi(\rd y).
$$
Hence, if  event $\cD$ is realized
\begin{eqnarray}
\label{eq4:AUX}
&&\Upsilon\big(X^{(n)}\big)=\prod_{m\in\cM}\frac{\gamma_{m,\mu}\big(X^{(n)}\big)}{\gamma_{m,\nu}\big(X^{(n)}\big)}=
\prod_{k=0}^{t-1}\bigg[\frac{\mT_k(\mu)}{\mT_k(\nu)}\bigg]^{\eta_k}.
\end{eqnarray}
Setting for $k=0,1,\ldots,t-1$
$$
u_{k}(\pi):= \sum_{j=0}^{t-k-1} \tfrac{(-1)^{j}}{j!}(D\widehat{n}_{0})^{j} e_{\pi}(j+k),\;\;\;
U_k(\pi):=
\sum_{j=t-k}^\infty \tfrac{(-1)^{j}}{j!} (D\widehat{n}_{0})^{j}e_{\pi}(j+k),
$$
we get in view of the Taylor expansion
$\mT_k(\pi)=u_{k}(\pi)+U_k(\pi)$. Moreover,
since $\mu\stackrel{t,t}{\backsim}\nu$,
\begin{equation}
\label{eq6:AUX}
u_{k}(\mu)=u_{k}(\nu),\;\;\;\forall k=0,1, \ldots, t-1.
\end{equation}
Next, in view of (\ref{eq9:AUX}),  $Dn<2\kappa_1<1$  so that   for sufficiently small $\kappa_1$
\begin{equation}
\label{eq5:AUX}
\big|U_{k}(\mc)\big|\leq (eDn)^{t-k} e_\pi(t), \quad k=0,1,\ldots,t-1.
\end{equation}
Also, by (\ref{eq9:AUX}), and by definition
 $k=0,1,\ldots,t-1$
\begin{equation}
\label{eq7:AUX}
\mT_k(\pi)\geq  e^{-Dn}e_{\pi}(k)\geq
e^{-2\kappa_1} e_\pi(t).
\end{equation}
We deduce from (\ref{eq6:AUX}), (\ref{eq5:AUX}) and (\ref{eq7:AUX}) that on the event  $\cD$
one has
\begin{eqnarray}
\label{eq8:AUX}
&&\frac{\mT_k(\mu)}{\mT_k(\nu)}\geq 1-2e^{2\kappa_1}(eDn)^{t-k}
\geq 1-4e^{2\kappa_1+1}\kappa_1 M^{(k-t)/t}.
\end{eqnarray}
To get the last inequality we  used (\ref{eq9:AUX}),  $Dn<2\kappa_1<1$, and took into account that $\kappa_1$ is small enough so that $2\kappa_1 e<1$.
Assuming additionally  that $\cE$ is realized
we get from (\ref{eq4:AUX}) and (\ref{eq8:AUX})
\begin{eqnarray*}
\label{eq11:AUX}
\Upsilon\big(X^{(n)}\big)&\geq&\prod_{k=0}^{t-1}
\Big[1- 4\kappa_1e^{2\kappa_1+1} M^{(k-t)/t}\Big]^{M^{(t-k)/t}}
\nonumber\\
&\geq&\Big[\inf_{z\geq 1}\Big(1-4\kappa_1e^{2\kappa_1+1}z^{-1}\Big)^{z}\Big]^{t}>1/2,
\end{eqnarray*}
provided that $\kappa_1$ is sufficiently small, and $t$ is independent of $n$.

Thus, we have proved that
\begin{equation*}
\label{eq12:AUX}
\cD\cap\cE\subseteq \big\{\Upsilon\big(X^{(n)}\big)\geq \tfrac{1}{2}\big\},
\end{equation*}
and, to complete the proof  it suffices to show that
\begin{equation}
\label{eq1200:AUX}
\bE_{\nu}\big[\bP_{f_{\bze}}\left\{\cD\cap\cE \right\}\big]\geq 1/3.
\end{equation}

2$^0$.  Let  $\bar{\cE}$ and $\bar{\cB}$ be  the events  complimentary to $\cE$ and $\cB$ respectively.
\par
(a). First we show that
\begin{equation}
\label{eq1201:AUX}
\bE_{\nu}\big[\bP_{f_{\bze}}\{\bar{\cE}\}\big]\leq e\kappa_1.
\end{equation}
By Markov's inequality
$$
\bP_{f_{\bze}}\left\{\bar{\cE} \right\}\leq \sum_{k=1}^{t-1}\bP_{f_{\bze}}\left\{\eta_k > M^{1-\frac{k}{t}}\right\}\leq \sum_{k=1}^{t-1} M^{\frac{k}{t}-1}\bE_{f_{\bze}}\left(\eta_k\right),
$$
where in the first inequality  we have used that $\eta_0\leq M$ by definition.
Noting that
$$
\bE_{f_{\bze}}\left(\eta_k\right)=\sum_{m\in\cM}\bP_{f_{\bze}}\{\widehat{n}_{m}=k\}=
\sum_{m\in\cM}\bP_{f_{\bze}}\bigg\{\sum_{i=1}^n\mathrm{1}_{\Pi_m}(X_i)=k\bigg\},
$$
we obtain
$$
\bE_{f_{\bze}}(\eta_k)=\sum_{m\in\cM}\tbinom{n}{k}
\Big[\int_{\Pi_m}f_{\bze}(x)\rd x\Big]^{k}\Big[1-\int_{\Pi_m}
f_{\bze}(x)\rd x\Big]^{n-k}\leq \frac{M(nA\bma)^k}{k!}.
$$
In the last inequality we took into account that  $\bze_m\leq 1$ for all $m\in\cM$.
Using the condition  (\ref{eq:condition-main}) we obtain
$\bE_{f_{\bze}}\left(\eta_k\right)=\kappa_1M^{(t-k)/t}(k!)^{-1}$,
and (\ref{eq1201:AUX}) follows.
\par
(b). Now let us prove that
\begin{equation}
\label{eq1202:AUX}
\bE_{\nu}\Big[\bP_{f_{\bze}}\left\{\bar{\cD} \right\}\Big]\leq (e\kappa_1)^{t}.
\end{equation}
Indeed,
$$
\bP_{f_{\bze}}\{\bar{\cD}\}\leq \sum_{m\in\cM}\bP_{f_{\bze}}\left\{\widehat{n}_{m}\geq t\right\}=
\sum_{m\in\cM}\bP_{f_{\bze}}\Big\{\sum_{i=1}^n\mathrm{1}_{\Pi_{m}}(X_i)\geq t\Big\}.
$$
Noting that under $f_\bze$  random variable $\widehat{n}_m$ is binomial with parameters $n$ and
$$
p_n:=\bP_{f_{\bze}}\big(X_i\in\Pi_{m}\big)=\int_{\Pi_m}f_{\bze}(y)\rd y=\mA\bma\bze_{m}\leq \mA\bma,
$$
$np_n\leq \kappa_1M^{-1/t}\leq t$ in view of the  condition (\ref{eq:condition-main}), and
applying Lemma~\ref{lem:lb-deviation} given in the proof of
Corollary~\ref{cor:lower-bound-genera2} we obtain for any $m\in\cM$
\begin{equation*}
\bP_{f_{\bze}}\big(\widehat{n}_{m}\geq t\big)\leq
\big(e n\mA\bma t^{-1}\big)^{t}\leq (e \kappa_1)^{t}M^{-1}.
\end{equation*}
This implies  (\ref{eq1202:AUX}), and it follows from  (\ref{eq1201:AUX}) and (\ref{eq1202:AUX})
that
$$
\bE_{\nu}\big[
\bP_{f_{\bze}}\{\cD\cap\cE \}\big]\geq 1-\big[e\kappa_1+(e\kappa_1)^{t}\big]>1/3,
$$
provided that $\kappa_1$ sufficiently small.
The corollary statement follows now from (\ref{eq1200:AUX}) and Theorem~\ref{th:lower-bound-general}.
\epr
\section{Appendix}
\subsection{ Proof of Proposition \ref{prop1}}
Fix $s\in \bN^*$ and let
$M_s$
denote the Hilbert matrix, that is $M_s=\big\{(i+j+1)^{-1}\big\}_{i,j=0,\ldots, s}$. It is well known that
$M_s$ is invertible for all $s\in\bN^*$.
\par
Let $\Ba=(a_0,\ldots,a_{s})\in\bR^{s+1}$ be such that
$$
\int_{0}^1P_{s,\Ba}(x)x^k\rd x=\delta_{k,t},\quad \forall k=0,\ldots s;
$$
here  $P_{s,\Ba}(x)=\sum_{j=0}^s a_j x^j$,   and $\delta_{k,t}$ is the Kronecker symbol.
If $\Be_t=(0, \ldots, 0, 1, 0,\ldots 0)\in \bR^{s+1}$ is the
$t$--th unit vector of the canonical basis of $\bR^{s+1}$ then
%
$$
\int_{0}^1P_{s,\Ba}(x)x^k\rd x=\delta_{k,t},\quad \forall k=0,\ldots s\quad \Leftrightarrow\quad M_{s}\Ba=
\Be_t,\;\;\;\Ba=M^{-1}_s \Be_t.
$$
%
\par Put for brevity $P=P_{s,\Ba}$, and  note that
\begin{equation}
\label{eq1}
0< \|P\|^2_{\bL_1(0,1)}\leq \|P\|^2_{\bL_2(0,1)}=\Ba^TM_{s}\Ba
=\Ba^T\Be_{t}=\big[M^{-1}_s]_{t,t}.
\end{equation}
Define
$P^+(x):=\max \{P(x),0\}$,  $P^-(x):=\max \{-P(x),0\}$
and remark that since  $P(x)=P^+(x)-P^+(x)$ for all $x\in[0,1]$
$$
0=\int_0^1 P(x)\rd x=\int_0^1 P^+(x)\rd x-\int_0^1 P^-(x)\rd x.
$$
Moreover $|P(x)|=P^+(x)+P^-(x)$ for any $x\in[0,1]$ and, therefore,
$$
\int_0^1 P^+(x)\rd x=\int_0^1 P^-(x)\rd x=\tfrac{1}{2}\|P\|_{\bL_1(0,1)}.
$$
Setting
$$
\mu^\prime(\rd x)=\frac{2P^+(x)\mathrm{1}_{[0,1]}(x)}{\|P\|_{\bL_1(0,1)}}\rd x,\quad
\nu^\prime(\rd x)=\frac{2P^-(x)\mathrm{1}_{[0,1]}(x)}{\|P\|_{\bL_1(0,1)}}\rd x
$$
we can assert that $\mu^\prime,\nu^\prime\in\mP[0,1]$.
Note also that in view of (\ref{eq1}) for any $k=1,\ldots,s$
\begin{equation}
\label{eq1-new}
e_{\mu^\prime}(k)-e_{\nu^\prime}(k)=2\|P\|^{-1}_{\bL_1(0,1)}\int_0^1x^k P(x)\rd x
=
2\|P\|^{-1}_{\bL_1(0,1)}\delta_{k,t}.
\end{equation}
Finally, let
$$
\mu=\tfrac{1}{2}\mu^\prime(\rd x)+\tfrac{1}{2}\delta_1(\rd x),\quad \nu=\frac{1}{2}\nu^\prime(\rd x)+\tfrac{1}{2}\delta_1(\rd x),
$$
where $\delta_t$, $t\in\bR$, denotes the  Dirac mass  at $t$.
It is obvious that (\ref{eq1-new}) is fulfilled for $\mu$ and $\nu$
with the constant $\|P\|^{-1}_{\bL_1(0,1)}$ if $k=t$. Additionally
all moments w.r.t. $\mu$ and $\nu$ are greater than $1/2$.
It remains to note that the elements of the inverse Hilbert matrix are known explicitly  and in particular
$
\big[M^{-1}_s]_{t,t}=C^{-2}_{s,t}.
$
This, together with (\ref{eq1}) completes the proof.
\epr
\subsection{\bf Proof of Proposition \ref{prop2}}
 1$^0$. The first step in the proof is to find a function $K:\bR\to\bR$ satisfying
\begin{eqnarray}
\label{eq2}
&&\int_0^1 K(x)x^k\rd x=0,\; k=0,\ldots,s, \quad \int_0^1 K(x)S(x)\rd x=\varpi_{s}(S).
\end{eqnarray}
We will seek $K$ in the following form:
$$
K(x)=P_{s,\Ba}(x)-bS(x),
$$
where $\Ba=(a_0,\ldots,a_s)\in\bR^{s+1}$,  $b\in\bR$  are the parameters to be chosen.
Let
$$
c_k=\int_{0}^1S(x)x^k\rd x, \quad k=0,\ldots,s,
$$
and let
$\Bc=\big(c_0,\ldots,c_{r}\big)\in\bR^{s+1}$.
Note that condition
$\int_0^1 K(x)x^k\rd x=0$, $k=0,\ldots,s$ is equivalent to   $M_s\Ba=b\Bc$, where as before  $M_s$ denotes
the Hilbert matrix.

Since $M_s$ is invertible we get
\begin{equation}\label{eq:Ba}
\Ba=bM_s^{-1}\Bc.
\end{equation}
With this choice the second condition in (\ref{eq2}) becomes
$$
b\Ba^T\Bc-b\|S\|_{\bL_2(0,1)}=b\Big[\Bc^TM_s^{-1}\Bc-\|S\|_{\bL_2(0,1)}\Big]=\varpi_s(S).
$$
It remains to note that
$$
\kappa_s(S):=\Bc^TM_s^{-1}\Bc-\|S\|_{\bL_2(0,1)}=-\inf_{\Bu\in\bR^{s+1}}\int_0^1\big|S(x)-P_{s,\Bu}(x)\big|^2\rd x,
$$
and $\kappa_s(S)<0$ because $\kappa_s(S)=0$ implies
that $S(x)$ coincides almost everywhere with a polynomial of degree $s$ which  contradicts to the proposition
assumption $\varpi_s(S)>0$. Choosing
$b=\varpi_s(S)/\kappa_s(S)$ we conclude that
$K(x)=P_{s, \Ba}(x)-b S(x)$ satisfies (\ref{eq2}).
\par\medskip
\noindent 2$^0$.
Let $\mL$ be the linear subspace of $\bC(0,1)$  spanned by
functions $\big\{S(x),1,x,x^2, \ldots, x^s\big\}$.
Let  $K(x)=P_{s,\Ba}(x)- [\varpi_s(S)/\kappa_s(S)]S(x)$, where $\Ba$ is defined in
(\ref{eq:Ba}). Define
$$
\Lambda(\ell)=\int_0^1 K(x)\ell(x)\rd x,\quad \ell\in\mL;
$$
$\Lambda$ is a linear continuous functional on
$\mL$, and its norm is
$\|\Lambda\|_{\mL}:=\sup_{\ell\in\mL}\{|\Lambda(\ell)|/\|\ell\|_\infty\}$.
Our goal now is to prove that
\begin{equation}
\label{eq3}
\|\Lambda\|_{\mL}=1.
\end{equation}
For any $\e>0$ let $P^{(\e)}$ be a polynomial of the degree $s$ such that
$$
\sup_{x\in[0,1]}\big|S(x)-P^{(\e)}(x)\big|\leq \varpi_s(S) (1+\e).
$$
Putting  $\ell_\e(x)=S(x)-P^{(\e)}(x)$ and noting that $\ell_\e\in\mL$ for any
$\e>0$ we have by (\ref{eq2})
\begin{eqnarray*}
\|\Lambda\|_{\mL}&\geq& \big|\Lambda(\ell_\e)\big|\,\|\ell_\e\|_\infty^{-1}= |\ell_\e\|^{-1}_\infty
\bigg|\int_0^1K(x)\ell_\e(x)\rd x\bigg|
\\
&=&\|\ell_\e\|^{-1}_\infty\bigg|\int_0^1 K(x)S(x)\rd x\bigg| =\varpi_s(S)
\big\|\ell_\e\|^{-1}_\infty\geq (1+\e)^{-1}.
\end{eqnarray*}
Since $\e$ can be chosen  arbitrary small we prove that necessarily $\|\Lambda\|_{\mL}\geq 1$.
\par
Assume that $\|\Lambda\|_{\mL}> 1$; then there exists $\ell_0\in\mL$ such that
$
|\Lambda(\ell_0)|> \|\ell_0\|_\infty.
$
Next we note that
$
\Lambda(\ell)=0
$
if and only if $\ell$ is a polynomial. Indeed, by definition of $\mL$ any $\ell\in\mL$
is represented as follows
$$
\ell(x)=t_{s+1}S(x)+P_{s,\Bt}(x),\quad \Bt=(t_0,\ldots,t_s)\in\bR^{s+1},\; t_{s+1}\in\bR.
$$
In view of (\ref{eq2}),
$\Lambda(\ell)= t_{s+1}\varpi_s(S)$  and, therefore,
$
\Lambda(\ell)=0
$ if and only if $t_{s+1}=0$.
Now we define
\[
l(x):=S(x)-\frac{\varpi_s(S)\ell_0(x)}{\Lambda(\ell_0)},
\]
and note that
$\Lambda(l)=\Lambda(S)-\varpi_s(S)=0$ which means that $l$ is a polynomial.
On the other hand,
$$
\sup_{x\in[0,1]}|S(x)-l(x)|=\|\varpi_s(S) \Lambda^{-1}(\ell_0)\ell_0\|_\infty=
\varpi_s(S) \big|\Lambda(\ell_0)\big|^{-1}\|\ell_0\|_\infty<\varpi_s(S),
$$
which  is impossible by definition of $\varpi_s(S)$. Thus, (\ref{eq3}) is established.
\par\medskip
\noindent 3$^0$. By the Hahn-Banach extension theorem there exists a linear continuous  functional $\Lambda^*$ on $\bC(0,1)$
satisfying
\begin{equation}
\label{eq4}
\Lambda^*(\ell)=\Lambda(\ell),\;\; \forall \ell\in\mL,\quad \|\Lambda^*\|_{\bC(0,1)}= \|\Lambda\|_{\mL}.
\end{equation}
The Riesz representation theorem implies existence of  a unique signed measure $\lambda$  such that
\begin{equation}
\label{eq5}
\Lambda^*(u)=\int_{0}^1u (x)\lambda(\rd x),\quad\forall u\in\bC(0,1),\quad |\lambda|\big({[0,1]}\big)=
\|\Lambda^*\|_{\bC(0,1)}.
\end{equation}
Moreover, in view of the Jordan decomposition theorem,
$\lambda$ can be represented uniquely as $\lambda=\lambda^+-\lambda^-$, where
$\lambda^+$, $\lambda^-$ are positive measures, and $|\lambda|=\lambda^+ + \lambda^-$.
Therefore we obtain from
(\ref{eq2}), (\ref{eq4}) and (\ref{eq5})
$$
\lambda([0,1])=\lambda^{+}\big({[0,1]}\big)-\lambda^{-}\big({[0,1]}\big)=\Lambda^*\big(1)=\Lambda(1)=0.
$$
In addition, (\ref{eq3}), (\ref{eq4}) and (\ref{eq5}) imply that
$
\lambda^{+}\big({[0,1]}\big)+\lambda^{-}\big({[0,1]}\big)=1
$
and, therefore,
\begin{equation}
\label{eq5-new}
\lambda^{+}\big({[0,1]}\big)=\lambda^{-}\big({[0,1]}\big)=\tfrac{1}{2}.
\end{equation}
Note that in view of  (\ref{eq2}), (\ref{eq4}) and (\ref{eq5}) for any $k=0,\ldots, s$ one has
$$
\int_0^1x^k\lambda^{+}(\rd x)-\int_0^1x^k\lambda^{-}(\rd x)=\int_0^1x^k\lambda(\rd x)=
\Lambda^*\big(x\mapsto x^k\big)=\Lambda\big(x\mapsto x^k\big)=0,
$$
and
$$
\int_0^1S(x)\lambda^{+}(\rd x)-\int_0^1S(x)\lambda^{-}(\rd x)=
\int_0^1S(x)\lambda(\rd x)=\Lambda^*(S)=\Lambda(S)=\varpi_s(S).
$$
Finally, in view of (\ref{eq5-new}), $2\lambda^{+}$ and
$2\lambda^{-}$ are probability measures on $[0,1]$; therefore
letting
$$
\mu(\rd x):=\lambda^+(\rd x)+ \tfrac{1}{2}\delta_1(\rd x),\quad \nu(\rd x):=\lambda^-(\rd x)+\frac{1}{2}\delta_1(\rd x)
$$
we complete the proof.
\epr

\subsection{Proof of Lemma \ref{lem:verification-ass-on-family}}

First of all we remark that the construction of the set of functions
$\big\{f_{w}, w\in \{0,1\}^{M}\big\}$ almost coincides with the construction proposed in \cite{gl14},
in the proof of Theorem~3.
Thus, denoting $F_{w}=\sum_{m\in\cM}w_{m}\Lambda_{m}$
and repeating computations done in the
cited paper we can verify that the assumption
\begin{eqnarray}
\label{eq8:proof-th:lower-bound-deconvolution}
&&\mA\sigma_l^{-\beta_l} \|\Lambda\|_{r_l}
\big[\varrho_{w}(r_l)\bma\big]^{1/r_l}\leq C_1 L_l,\quad\forall l=1,\ldots, d,
\end{eqnarray}
together with $\vec{\sigma}\in (0,1]^d$ guarantees $F_{w}\in\bN_{\vec{r},d}(\vec{\beta},\tfrac{1}{2}\vec{L})$, $w\in[0,1]^{M}$.
It is important to realize that the only conditions used in the
proof of (\ref{eq8:proof-th:lower-bound-deconvolution}) is (\ref{eq:Lambda}) and $\Lambda\in\bC^{\infty}(\bR^d)$ which
are the same as  in \cite{gl14}.  Set for any $z\geq 0$
\begin{eqnarray*}
\cW_{\mc,z}:=\Big\{w\in[0,1]^{M}:\; \big|\varrho_{w}(z)- M
e_{\mc}(z)\big|\leq \sqrt{\upsilon Me_\pi(2z)}\Big\}.
\end{eqnarray*}
First we note that if for some $z\geq 1$ the event $\big\{\bze\in\cW_{\mc,z}\big\}$ is realized, then
\begin{equation}
\label{eq1:proof-prop3}
\varrho_{\bze}(z)\leq M e_\mc(z)+\sqrt{\upsilon Me_\mc(2z)}=e_\mc(z)[M +\sqrt{2\upsilon M}\big]\leq 2Me_\mc(z)\leq 2Me^*(z),
\end{equation}
since $e_{\mc}(z)\leq 1$, $M\geq 36\upsilon$ in view of  (\ref{eq3**-lbf}). Also we have used (\ref{eq:e>1/2}).

\noindent Next, we deduce from  (\ref{eq1:proof-prop3})  that if  $\big\{\bze\in \cW_{\mc,z}\big\}$, $z\in[1,\infty)$ is realized then
\begin{equation*}
\big\|F_{\bze}\big\|_z= \|\Lambda\|_{z}\big(\bma\varrho_{\bze}(z))^{\frac{1}{z}}
\leq 2 \|\Lambda\|_{z}\big(\bma M e^*(z)\big)^{\frac{1}{z}}.
\end{equation*}
It yields in particular, $\big\{\bze\in W_{\mc,q}\big\}$ is realized
$$
\big\|f_{\bze}\big\|_q=\big[1-A\bma\varrho_{\bze}(1)\big]\|f_{0}\|_q+ A\big\|F_{\bze}\big\|_q\leq 2^{-1}Q+2\mA\|\Lambda\|_{q}\big(\bma M e^*(q)\big)^{\frac{1}{q}}\leq Q.
$$
To get the last inequality we have used (\ref{eq3*-lbf}) and that $f_{0}\in\bB_q(Q/2)$ in view of the second assertion of Lemma \ref{lem:density-f_0,N}.
Thus, $f_{\bze}\in\bB_q(Q)$.

\noindent At last, if $\big\{\bze\in \cap_{l=1}^d\cW_{\mc,r_l}\big\}$ is realized we have in view of condition (\ref{eq4-lbf}) and (\ref{eq1:proof-prop3})
$$
\mA\|\Lambda\|_{r_l}\big(\bma\varrho_{\bze}(r_l)\big)^{\frac{1}{r_l}}\leq 2\mA\|\Lambda\|_{r_l}\big(\bma M e^*(r_l)\big)^{\frac{1}{r_l}}\leq 2c_1 L_l
$$
for any $r_l\neq\infty$. Additionally, for all $l\in\{1,\ldots, d\}$ such that $r_l=\infty$ the latter inequality obviously hold for all realizations of $\bze$.

We  assert in view of (\ref{eq8:proof-th:lower-bound-deconvolution}) with $C_1=2c_1$  that $F_{\bze}\in \bN_{\vec{r},d}(\vec{\beta},2^{-1}\vec{L})$ and, therefore,
$f_{\bze}\in \bN_{\vec{r},d}(\vec{\beta},\vec{L})$ since $f_{\bze}=\big[1-A\bma\varrho_{\bze}(1)\big]f_{0}+F_{\bze}$ and $f_{0}\in \bN_{\vec{r},d}(\vec{\beta},2^{-1}\vec{L})$
in view of the second assertion of Lemma \ref{lem:density-f_0,N}. Thus, we have shown that
\begin{eqnarray*}
&&\bigcap_{z\in\{r_1,\ldots,r_d,q\}}\big\{\bze\in  \cW_{\mc,z}\big\}\subset \left\{f_{\bze}\in\bN_{\vec{r},d}(\vec{\beta},\vec{L})\cap\bB_q(Q)\right\}.
\end{eqnarray*}
 Hence
\begin{eqnarray*}
&&\bP_{\mc}\big\{f_{\bze}\in\bN_{\vec{r},d}(\vec{\beta},\vec{L})\cap\bB_q(Q)\big\}\geq 1-\sum_{z\in\{r_1,\ldots,r_d,q\}}\bP_{\mc}\big\{\bze\notin  \cW_{\mc,z}\big\}
\\
&&\geq 1-(d+1)\upsilon^{-1} =1-(64)^{-1}.
\end{eqnarray*}
Lemma is proved.
\epr

\bibliographystyle{agsm}

\end{document}